\documentclass{article}
\usepackage{mathrsfs}
\usepackage{amsfonts}
\usepackage{multirow}
\usepackage[dvips]{graphicx}

\usepackage{amsmath,amssymb}
\usepackage{mathrsfs,color}
\newenvironment{proof}[1][Proof]{\noindent\textbf{#1.} }{\ \rule{0.5em}{0.5em}}

\numberwithin{figure}{section} \numberwithin{equation}{section}
\makeatletter 
\newcommand{\red}[1]{{\color{red}{#1}}}

\begin{document}
\date{}
\author{Hengfei Ding$^{1}$,
\;\;Changpin
Li$^{2}$
\\
 \small \textit{1. School of Mathematics and Statistics, Tianshui
Normal University, Tianshui 741001, China}\\
\small \textit{2. Department of Mathematics, Shanghai University,
Shanghai 200444, China }
} \vspace{0.2cm}
\title{ High-order Compact Difference Schemes for the Modified Anomalous Subdiffusion
Equation\thanks{~The work was partially supported by the National
Natural Science Foundation of China under Grant No. 11372170,
Key Program of Shanghai Municipal Education Commission
under Grant No.  12ZZ084, the grant of ``The First-class
Discipline of Universities in Shanghai", the Scientific Research Program
 for Young Teachers of Tianshui Normal University under Grant No. TSA1405, and Tianshui
Normal University Key Construction Subject Project (Big data processing in dynamic image)}.
}\maketitle \vspace{0.2 cm}
 \hrulefill

\begin{abstract}
In this paper, two kinds of high-order compact finite difference
schemes for second-order derivative are developed. Then a second-order numerical scheme for a Riemann-Liouvile
derivative is established based on a fractional  centered difference operator. We apply
 these methods to a fractional anomalous subdiffusion
equation to construct two kinds of novel numerical schemes. The solvability,
stability and convergence analysis of these difference schemes are
studied by using Fourier method. The convergence orders of
 these numerical schemes are $\mathcal {O}(\tau^2+h^6)$ and $\mathcal
{O}(\tau^2+h^8)$, respectively.
 Finally, numerical experiments are displayed which are in line with the theoretical analysis.
\\\vspace{0.1 cm}
\vspace{0.1 cm}\\
\end{abstract}
\hrulefill

\section{Introduction}
\quad\ The phenomenological diffusion equation can be derived by the following
Fick's first law \cite{CM} (which describes steady-state diffusion),
$$
\begin{array}{lll} \displaystyle
\mathcal {J}=-\kappa\nabla u.
\end{array}\eqno(1)
$$
Combing the following conservation law of energy
$$
\begin{array}{lll} \displaystyle
\frac{\partial u}{\partial t}=-\nabla\cdot\mathcal {J},
\end{array}\eqno(2)
$$
one can obtain the diffusion equation below (also known as Fick's second law or the heat equation)
$$
\begin{array}{lll} \displaystyle
\frac{\partial u}{\partial t}=\nabla\cdot(\kappa\nabla u).
\end{array}
$$
This equation well characterizes the classic diffusion phenomenon \cite{CR}.

However, if the diffusion is abnormal,
that is to say, it follows non-Gaussian statistics or can be
interpreted as the L\'{e}vy stable densities, then
the above equation can not well describe such anomalous diffusion. Generally
speaking, the fractional differential equations can well describe and
model these anomalous diffusion phenomena \cite{SCK}.
The corresponding fractional Fick's law has been proposed \cite{SK}.
$$
\begin{array}{lll} \displaystyle
\mathcal {J}_{\mathcal {A}}=-\mathcal
{A}\cdot\,_{RL}D_{0,t}^{1-\alpha}\nabla u,\;\alpha\in(0,1).
\end{array}
$$
Combination of this equation with equation (2) gives
$$
\begin{array}{lll} \displaystyle
\frac{\partial u}{\partial t}=\mathcal
{A}\cdot\,_{RL}D_{0,t}^{1-\alpha}\frac{\partial^2 u}{\partial x^2},
\end{array}
$$
where $0<\alpha<1$, and $\mathcal {A}>0$ is the anomalous
diffusion coefficient. If $\alpha=1$, it is just the normal diffusion equation. Here $\,_{RL}D_{0,t}^{1-\alpha}$ is the
Riemann-Liouville operator, which is defined as follows:
$$
\begin{array}{ll}
\displaystyle
\,_{RL}D_{0,t}^{1-\alpha}u(x,t)=\frac{1}{\Gamma(\alpha)}\frac{\partial}{\partial
t}\int_{0}^{t}\frac{u(x,s)}{(t-s)^{1-\alpha}}ds,
\end{array}
$$
where $\Gamma(\cdot)$ is the Gamma function.

Recently, a modified fractional Fick's law has been used to describe processes that become less anomalous
as time progresses by the inclusion of a secondary fractional time derivative acting on a diffusion operator \cite{CGSY},
$$
\begin{array}{lll} \displaystyle
\mathcal {J}_{\mathcal {A,B}}=-\left(\mathcal
{A}\cdot\,_{RL}D_{0,t}^{1-\alpha}+\mathcal
{B}\cdot\,_{RL}D_{0,t}^{1-\beta}\right) \nabla u,
\end{array}
$$
where $0<\alpha<1,0<\beta<1$, and $\mathcal {A}>0,\mathcal {B}>0$ are the anomalous
diffusion coefficients.
Thus, the modified fractional
anomalous diffusion equation is obtained \cite{L},
$$
\begin{array}{lll} \displaystyle
 \frac{\partial{{}u(x,t)}}{\partial{t}}
 =\left(\mathcal {A}\cdot\,_{RL}D_{0,t}^{1-\alpha}+
\mathcal {B}\cdot\,_{RL}D_{0,t}^{1-\beta}
 \right)\left[
 \frac{\partial^2 u(x,t)}{\partial{x^2}}\right].
\end{array}
$$

 Till now, various kinds of anomalous diffusion equations have been studied numerically,
see $\cite{WWS,CLTA,CLB,ZLL,YA,ZTD,MT,HTVY,CJLT,WZ,AAA}$
and many the references cited therein. However, it seems that only a few numerical studies are available for the
two-term subdiffusions of the above form.

In the present paper, we aim to study the
following modified anomalous diffusion equation with a source term
$$
\begin{array}{lll} \displaystyle
 \frac{\partial{{}u(x,t)}}{\partial{t}}
 =\left(\mathcal {A}\cdot\,_{RL}D_{0,t}^{1-\alpha}+
\mathcal {B}\cdot\,_{RL}D_{0,t}^{1-\beta}
 \right)\left[
 \frac{\partial^2 u(x,t)}{\partial{x^2}}\right]
  +f(x,t),\;\;\;0<x<L,\;\;\;0< t\leq T$$,
\end{array}\eqno(3)
$$
subject to the initial and Dirichlet boundary value conditions
$$
\begin{array}{ll}
u(x,0)=\phi(x),\;\;0<x<L,
\end{array}
$$
$$
\begin{array}{ll}
u(0,t)=\varphi_{_1}(t),\;\;0\leq t\leq T,
\end{array}
$$
$$
\begin{array}{ll}
u(L,t)=\varphi_{_2}(t),\;\;0\leq t\leq T,
\end{array}
$$
where $f(x,t)$, $\phi(x)$, $\varphi_{_1}(t)$ and $\varphi_{_2}(t)$
are suitably
 smooth.

Jiang and Chen
proposed a collocation method based on reproducing kernels to
solve a modified anomalous subdiffusion equation (3)
with a linear source term on a finite domain \cite{JC}. In
\cite{LYB}, Liu et al., constructed a conditionally
stable difference scheme for equation (3) with a nonlinear source term, and
they proved that the convergence order is $\mathcal
{O}(\tau+h^2)$ by the energy method. In \cite{MAD}, Mohebbi et al.
considered an unconditionally stable difference scheme of order
$\mathcal {O}(\tau+h^4)$. Wang and
Vong \cite{WV} presented a compact method for the numerical
simulation of the modified anomalous subdiffusion equation (3), and
they achieved the convergence order $\mathcal
{O}(\tau^2+h^4)$. The aim of this paper is to propose
much higher order numerical methods for equation (3). We
construct two kinds of high-order compact difference schemes and provide a detailed study of
the stability and convergence of the proposed methods by using the
Fourier method. We demonstrate that the convergence orders are $\mathcal
{O}(\tau^2+h^6)$ and $\mathcal {O}(\tau^2+h^8)$, respectively. One of advantages of compact difference schemes are that they can produce highly accurate numerical solutions but involves the less number of grid points. Thus, compact schemes result in matrices that have smaller band-width compared with non-compact schemes. For example, a sixth-order finite difference scheme involves seven grid points, while sixth-order compact difference scheme only needs five grid points. Another additional advantage of the compact high order methods is that the methods described here leads to diagonal linear systems, thus allowing the use of fast diagonal solvers. Different from the typical differential equations, even if we use the lower order methods for solving the fractional differential equations, we still need more calculations and strong spaces. If we use the higher order methods for fractional differential equations, the calculations and memory capacities can not remarkably increase. In this sense, the higher order numerical methods for fractional calculus and fractional differential equations attract more and more interest.

The rest of this article is organized as follows. In Section 2, we
firstly develop a sixth-order and an eight-order difference scheme
for second-order derivative, next a second-order numerical scheme
for the Riemann-Liouville derivative is proposed. Applications of
these methods to equation (3) give
 two effective finite difference schemes. The solvability,
stability and convergence of the numerical methods are discussed in
Sections 3, 4 and 5, respectively. The numerical experiments are performed
for equation (3) with the methods developed in this paper are given in Section 6,
which support the theoretical analysis. Finally, concluding
remarks are drawn in the last section.

\section{Numerical Schemes}
\quad\;\ Let
$t_k=k\tau~(k=0,1,\cdots,N)$
 and
$x_j=jh~(j=0,1,\cdots,M),$
 where the grid sizes in time and space are
defined by $\tau = T/N$ and $h = L/M$,  respectively.

Define the following centered difference operator as
$$\begin{array}{lll}
\displaystyle \delta_x
u(x_j,t_k)=u(x_{j+\frac{1}{2}},t_k)-u(x_{j-\frac{1}{2}},t_k),
\end{array}$$
then we have
$$\begin{array}{lll}
\displaystyle \delta_x^2
u(x_j,t_k)=u(x_{j+1},t_k)-2u(x_j,t_k)+u(x_{j-1},t_k).
\end{array}
$$

 It is well known that a
second-order approximation for the derivative $\frac{\partial^2
u(x_j,t_k)}{\partial{x^2}}$ is given by the following second-order centered
difference scheme
$$\begin{array}{rrr}
\displaystyle \frac{\partial^2
u(x_j,t_k)}{\partial{x^2}}=\frac{\delta_x^2
u(x_j,t_k)}{h^2}+\mathcal {O}(h^2).
\end{array}$$

A fourth-order compact difference scheme has also been
constructed \cite{C},
$$\begin{array}{rrr}
\displaystyle \frac{\partial^2
u(x_j,t_k)}{\partial{x^2}}=\frac{1}{h^2}\left(1+\frac{1}{12}\delta_x^2\right)^{-1}

\delta_x^2
u(x_j,t_k)+\mathcal
{O}(h^4).
\end{array}$$

Next, we develop two high-order compact difference schemes for the
second-order spatial derivative by the following lemma.

 {\it{\bf Lemma 1.
} Define the following two operators:
$$\begin{array}{lll}
\displaystyle\mathscr{L}_1=:\frac{1}{h^2}\left(1-\frac{1}{90}\delta_x^4\right)^{-1}
\delta_x^2\left(1-\frac{1}{12}\delta_x^2\right),
\end{array}$$
and
$$\begin{array}{lll}
\displaystyle \mathscr{L}_2=:\frac{1}{h^2}\left(1+\frac{1}{560}\delta_x^6\right)^{-1}
\delta_x^2\left(1-\frac{1}{12}\delta_x^2+\frac{1}{90}\delta_x^4\right),
\end{array}$$
then
$$\begin{array}{lll}
\displaystyle\frac{\partial^2
u(x_j,t_k)}{\partial{x^2}}=\mathscr{L}_1 u(x_j,t_k)+\mathcal
{O}(h^6)
\end{array}\eqno(4)
$$
and
$$\begin{array}{lll}
\displaystyle\frac{\partial^2
u(x_j,t_k)}{\partial{x^2}}=\mathscr{L}_2 u(x_j,t_k)+\mathcal
{O}(h^{8})
\end{array}\eqno(5)
$$
hold.}

\begin{proof}
In view of the following approximation scheme \cite{A}
$$\begin{array}{lll}
\displaystyle \frac{\partial^2
u(x_j,t_k)}{\partial{x^2}}&=&\displaystyle\left[\frac{2}{h}\sinh^{-1}\left(\frac{\delta_x}{2}\right)\right]^2
u(x_j,t_k)\vspace{0.2
cm}\\&=&\displaystyle\frac{1}{h^2}\left[\delta_x-\frac{1}{24}\delta_x^3+\frac{3}{640}\delta_x^5
-\frac{5}{7168}\delta_x^7+\cdots \right]^2u(x_j,t_k)\vspace{0.2
cm}\\&=&\displaystyle\frac{1}{h^2}\left[\delta_x^2-\frac{1}{12}\delta_x^4+\frac{1}{90}\delta_x^6
-\frac{1}{560}\delta_x^8+\frac{1}{3150}\delta_x^{10}-\frac{1}{16632}\delta_x^{12}+\cdots
\right]u(x_j,t_k),
\end{array}$$
then one obtains
$$\begin{array}{lll}
&\displaystyle \frac{1}{h^2}\left(1-\frac{1}{90}\delta_x^4\right)^{-1}\delta_x^2\left(1-\frac{1}{12}\delta_x^2\right)
u(x_j,t_k)\vspace{0.3 cm}\\&=\displaystyle\frac{1}{h^2}
\left[\delta_x^2-\frac{1}{12}\delta_x^4+\frac{1}{90}\delta_x^6
-\frac{1}{1080}\delta_x^8+\cdots \right]u(x_j,t_k)\vspace{0.2
cm}\\&=\displaystyle\frac{\partial^2
u(x_j,t_k)}{\partial{x^2}}+\frac{13}{15120h^2}\delta_x^8u(x_j,t_k)+\mathcal
{O}(h^8)\vspace{0.2 cm}\\&=\displaystyle\frac{\partial^2
u(x_j,t_k)}{\partial{x^2}}+\frac{13h^6}{15120}\frac{\partial^8
u(x_j,t_k)}{\partial{x^8}}+\mathcal {O}(h^8)
\end{array}$$
and
$$\begin{array}{lll}
\displaystyle
&\displaystyle\frac{1}{h^2}\left(1+\frac{1}{560}\delta_x^6\right)^{-1}
\delta_x^2\left(1-\frac{1}{12}\delta_x^2+\frac{1}{90}\delta_x^4\right)
u(x_j,t_k)\vspace{0.3cm}\\&=\displaystyle\frac{1}{h^2}
\left[\delta_x^2-\frac{1}{12}\delta_x^4+\frac{1}{90}\delta_x^6
-\frac{1}{560}\delta_x^8+\frac{1}{6720}\delta_x^{10}+\cdots
\right]u(x_j,t_k)\vspace{0.2 cm}\\&=\displaystyle\frac{\partial^2
u(x_j,t_k)}{\partial{x^2}}+\frac{197}{28672000h^2}\delta_x^{10}u(x_j,t_k)+\mathcal
{O}(h^{10})\vspace{0.2 cm}\\&=\displaystyle\frac{\partial^2
u(x_j,t_k)}{\partial{x^2}}+\frac{197h^8}{28672000}\frac{\partial^{10}
u(x_j,t_k)}{\partial{x^{10}}}+\mathcal {O}(h^{10}).
\end{array}$$
That is,
$$\begin{array}{lll}
\displaystyle\frac{\partial^2
u(x_j,t_k)}{\partial{x^2}}=\frac{1}{h^2}\left(1-\frac{1}{90}\delta_x^4\right)^{-1}
\delta_x^2\left(1-\frac{1}{12}\delta_x^2\right)
u(x_j,t_k)+\mathcal
{O}(h^6),
\end{array}$$
and
$$\begin{array}{lll}
\displaystyle\frac{\partial^2 u(x_j,t_k)}{\partial{x^2}}=
\frac{1}{h^2}\left(1+\frac{1}{560}\delta_x^6\right)^{-1}
\delta_x^2\left(1-\frac{1}{12}\delta_x^2+\frac{1}{90}\delta_x^4\right)
u(x_j,t_k)+\mathcal
{O}(h^{8}).
\end{array}$$
This completes the proof.
\end{proof}

{\it{\bf Lemma 2 \cite{QSS}} For the suitably smooth function $u(x,t)$
with respect to $x$, arbitrary different numbers $p$, $q$ and $s$,
one has
$$\displaystyle u(x,t_s)=\frac{(t_{s}-t_q)u(x,t_p)+(t_p-t_s)u(x,t_{q})}{t_{p}-t_q}+\mathcal
{O}\left(\left|(t_p-t_s)(t_q-t_s) \right|\right).$$}

Next, we develop a second order numerical scheme for the Riemann-Liouvile
derivative at nongrid points $\left(x_j,t_{k+\frac{1}{2}}\right)$.

In \cite{TG}, Tuan and Gorenflo introduced the following
fractional central difference operator
$$\begin{array}{lll}
\displaystyle
\Delta_{c,\tau}^{\gamma}u(x,t)=\sum\limits_{\ell=0}^{\infty}\varpi_\ell^{(\gamma)}
u\left(x,t-\left(\ell-\frac{\gamma}{2}\right)\tau\right),
\end{array}
$$
and proved that
$$\begin{array}{lll}
\displaystyle
\,_{RL}D_{0,t}^{\gamma}u\left(x,t\right)=\frac{1}{\tau^{\gamma}}
\sum\limits_{\ell=0}^{\infty}\varpi_\ell^{(\gamma)}
u\left(x,t-\left(\ell-\frac{\gamma}{2}\right)\tau\right)+\mathcal
{O}(\tau^2),
\end{array}\eqno(6)
$$
where $ \varpi_{\ell}^{(\gamma)}=(-1)^{\ell} \left(\gamma\atop \ell
\right) $.

Accordingly, we obtain the following form at point
$\left(x_j,t_{k+\frac{1}{2}}\right)$ in view of equation (6),
$$\begin{array}{lll}
\displaystyle
\,_{RL}D_{0,t}^{\gamma}u\left(x_j,t_{k+\frac{1}{2}}\right)=\frac{1}{\tau^{\gamma}}
\sum\limits_{\ell=0}^{\infty}\varpi_\ell^{(\gamma)}
u\left(x_j,t_k-\left(\ell-\frac{\gamma+1}{2}\right)\tau\right)+\mathcal
{O}(\tau^2).
\end{array}\eqno(7)
$$

Letting $t_s=t_k-\left(\ell-\frac{\gamma+1}{2}\right)\tau$,
$t_p=t_k-\left(\ell-1\right)\tau$ and $t_q=t_k-\ell\tau$ gives the following second-order numerical formula by using
equation (7) and Lemma 2,
$$\begin{array}{lll}
\displaystyle
\,_{RL}D_{0,t}^{\gamma}u\left(x_j,t_{k+\frac{1}{2}}\right)=&\displaystyle\frac{1}{2\tau^{\gamma}}
\sum\limits_{\ell=0}^{\infty}\varpi_\ell^{(\gamma)}\left((1+\gamma)
u\left(x_j,t_k-\left(\ell-1\right)\tau\right)\right.\vspace{0.2
cm}\\&\displaystyle\left.+(1-\gamma)
u\left(x_j,t_k-\ell\tau\right)\right) +\mathcal {O}(\tau^2).
\end{array}\eqno(8)
$$
Now set
$$ \displaystyle\tilde{u}(x,t)=\left\{
  \begin{array}{cc} \vspace{0.3 cm}
   u(x,t),\;\;\;t\in[0,T],\\
   \displaystyle  0,\;\;\;t\notin[0,T],
    \end{array}
   \right.$$
then the numerical formula (8) becomes
$$\begin{array}{lll}
\displaystyle
\,_{RL}D_{0,t}^{\gamma}u\left(x_j,t_{k+\frac{1}{2}}\right)&=&\displaystyle\frac{1+\gamma}{2\tau^{\gamma}}
\sum\limits_{\ell=0}^{k+1}\varpi_\ell^{(\gamma)}
u\left(x_j,t_k-\left(\ell-1\right)\tau\right)\vspace{0.2
cm}\\&&\displaystyle+\frac{1-\gamma}{2\tau^{\gamma}}
\sum\limits_{\ell=0}^{k}\varpi_\ell^{(\gamma)}
u\left(x_j,t_k-\ell\tau\right) +\mathcal {O}(\tau^2) \vspace{0.2
cm}\\&=&\displaystyle \frac{1}{\tau^{\gamma}}
\sum\limits_{\ell=0}^{k+1}g_\ell^{(\gamma)}
u\left(x_j,t_k-\left(\ell-1\right)\tau\right)+\mathcal {O}(\tau^2) ,
\end{array}\eqno(9)
$$
where
$$\begin{array}{lll}
\displaystyle
g_0^{(\gamma)}=\frac{1+\gamma}{2}\varpi_0^{(\gamma)},\;\;\;
g_\ell^{(\gamma)}=\frac{1+\gamma}{2}\varpi_\ell^{(\gamma)}+\frac{1-\gamma}{2}\varpi_{\ell-1}^{(\gamma)},\;\;\ell\geq1.
\end{array}
$$

Applying the Crank-Nicolson method to equation
(3) yields
$$
\begin{array}{lll} \displaystyle
 \frac{u(x_j,t_{k+1})-u(x_j,t_k)}{\tau}
 =&\displaystyle\left(\mathcal {A}\cdot\,_{RL}D_{0,t}^{1-\alpha}+
\mathcal {B}\cdot\,_{RL}D_{0,t}^{1-\beta}
 \right)\left[
 \frac{\partial^2
 u\left(x_j,t_{k+\frac{1}{2}}\right)}{\partial{x^2}}\right]\vspace{0.2 cm}\\ \displaystyle
 & +f\left(x_j,t_{k+\frac{1}{2}}\right)+\mathcal {O}(\tau^2).
\end{array}\eqno(10)
$$
Setting
$$
\begin{array}{lll} \displaystyle
w\left(x_j,t_{k+\frac{1}{2}}\right)=
 \frac{\partial^2
 u\left(x_j,t_{k+\frac{1}{2}}\right)}{\partial{x^2}}
\end{array}\eqno(11)
$$
and
substituting (9) into (10) leads to
$$
\begin{array}{lll} \displaystyle
 \frac{u(x_j,t_{k+1})-u(x_j,t_k)}{\tau}
 =
 \frac{\mathcal {A}}{\tau^{1-\alpha}}
\sum\limits_{\ell=0}^{k+1}g_\ell^{(1-\alpha)}
w(x_j,t_{k+1-\ell})\\
\displaystyle + \frac{\mathcal {B}}{\tau^{1-\beta}}
\sum\limits_{\ell=0}^{k+1}g_\ell^{(1-\beta)} w(x_j,t_{k+1-\ell})
  +f\left(x_j,t_{k+\frac{1}{2}}\right)+\mathcal {O}(\tau^2).
\end{array}\eqno(12)
$$

Let $u_j^k$ be the approximation solution of $u(x_j,t_k)$. Noting
equation (11) and substituting (4) and (5) into (12) give the following two finite difference schemes
for equation (3):
$$
\begin{array}{lll} \displaystyle
 \frac{u_j^{k+1}-u_j^k}{\tau}
 = \frac{\mathcal {A}}{\tau^{1-\alpha}}\mathscr{L}_1
\sum\limits_{\ell=0}^{k+1}g_\ell^{(1-\alpha)}u_j^{k+1-\ell}
+\frac{\mathcal {B}}{\tau^{1-\beta}}\mathscr{L}_1
\sum\limits_{\ell=0}^{k+1}g_\ell^{(1-\beta)}u_j^{k+1-\ell}
+f_j^{k+\frac{1}{2}}, \vspace{0.2 cm}\\\hspace{6.7cm} \displaystyle
 \;\;\;
0\leq k\leq N-1,\;1\leq j\leq M-1,\vspace{0.2 cm}\\
\displaystyle u_j^0=\phi(x_j),\;\;0\leq j\leq M,\vspace{0.2 cm}\\
\displaystyle
u_0^k=\varphi_1(t_k),\;\;u_M^k=\varphi_2(t_k),\;\;0\leq k\leq N;
\end{array}\eqno(13)
$$
and
$$
\begin{array}{lll} \displaystyle
 \frac{u_j^{k+1}-u_j^k}{\tau}
 = \frac{\mathcal {A}}{\tau^{1-\alpha}}\mathscr{L}_2
\sum\limits_{\ell=0}^{k+1}g_\ell^{(1-\alpha)}u_j^{k+1-\ell}
+\frac{\mathcal {B}}{\tau^{1-\beta}}\mathscr{L}_2
\sum\limits_{\ell=0}^{k+1}g_\ell^{(1-\beta)}u_j^{k+1-\ell}+f_j^{k+\frac{1}{2}},
\vspace{0.2 cm}\\\hspace{6.7cm} \displaystyle
  \;\;\;
0\leq k\leq N-1,\;1\leq j\leq M-1,\vspace{0.2 cm}\\
\displaystyle u_j^0=\phi(x_j),\;\;0\leq j\leq M,\vspace{0.2 cm}\\
\displaystyle
u_0^k=\varphi_1(t_k),\;\;u_M^k=\varphi_2(t_k),\;\;0\leq k\leq N.
\end{array}\eqno(14)
$$

It is obvious that the local truncation
errors of difference schemes (13) and (14) are $R_j^k=\mathcal
{O}(\tau^2+h^6)$ and $\widetilde{R}_j^k=\mathcal {O}(\tau^2+h^8)$,
respectively.

\section{Solvability Analysis}
Denote
$$
\begin{array}{lll}
\displaystyle {
\textbf{U}}^0=\left(\phi(x_1),\phi(x_2),\cdots,\phi(x_{M-1})\right)^{T},
{\textbf{U}}^k=\left(u_1^k,u_2^k,\cdots,u_{M-1}^k\right)^{T},\;k=1,2,\cdots,N,
\end{array}
$$
and
$$
\begin{array}{lll}
\displaystyle {
\textbf{F}}^k=\left(f_1^{k+\frac{1}{2}},f_2^{k+\frac{1}{2}},\cdots,f_{M-1}^{k+\frac{1}{2}}\right)^{T},\;k=0,1,\cdots,N.
\end{array}
$$

Then we obtain the matrix form of difference scheme (13)
$$
\begin{array}{lll}
\displaystyle
 \left(A-g_{0}^{(\alpha,\beta)}B\right)
 { \textbf{U}}^{k+1}=
 \left(A+g_{1}^{(\alpha,\beta)}B\right){
 \textbf{U}}^{k}+\sum\limits_{\ell=2}^{k+1}
g_{\ell}^{(\alpha,\beta)}B{
 \textbf{U}}^{k+1-\ell}+\tau
 A{\textbf{F}}^{k}+C_k,\vspace{0.2 cm}\\ \hspace{5cm}j=1,2,\ldots,M-1,
\;\;k=0,1,\cdots,N-1,
\end{array}\eqno(15)
$$
where $\displaystyle\mu_{\alpha}=\frac{\tau^{\alpha}}{h^2}\mathcal
{A},\;\;\mu_{\beta}=\frac{\tau^{\beta}}{h^2}\mathcal {B}$,
$g_{\ell}^{(\alpha,\beta)}=
\mu_{\alpha}g_{\ell}^{(1-\alpha)}
+\mu_{\beta}g_{\ell}^{(1-\beta)}\;(\ell=0,1,\ldots,k+1)$,
matrices $A, B, C_k$ are given in the Appendix I.

Similarly, the matrix form of the difference scheme (14) is given by
$$
\begin{array}{lll}
\displaystyle
 \left(\widetilde{A}-g_{0}^{(\alpha,\beta)}\widetilde{B}\right)
 { \textbf{U}}^{k+1}=
 \left(\widetilde{A}+g_{1}^{(\alpha,\beta)}\widetilde{B}\right){
 \textbf{U}}^{k}+\sum\limits_{\ell=2}^{k+1}
g_{\ell}^{(\alpha,\beta)}\widetilde{B}{
 \textbf{U}}^{k+1-\ell}+\tau A{\textbf{F}}^{k}+\widetilde{C}_k,\vspace{0.2 cm}\\ \hspace{5cm}j=1,2,\ldots,M-1,
\;\;k=0,1,\cdots,N-1,
\end{array}\eqno(16)
$$
where matrices $\widetilde{A}, \widetilde{B}, \widetilde{C_k}$ are also given in the Appendix I.

{{\it{\bf Remark 1.} In difference schemes (13) and (14), there
are some points $u_{-2}^k$, $u_{-1}^k$, $u_{M+1}^k$
and $u_{M+2}^k$ outside of the interval $[0,L]$, denoted as ghost-points, that are generally approximated using extrapolation
 formulas}, see Appendix II for more details.}

 {\it{\bf Lemma 3 \cite{GS}}. A circulant matrix $S$ is a Toeplitz matrix in
the form
$$ \displaystyle
S= \left(
  \begin{array}{cccccc}
   s_1 & s_2&  s_3& & \cdots& s_{M-1} \vspace{0.2 cm}\\
   s_{M-1} &  s_1 &s_2 & s_3& &\vdots \vspace{0.2 cm}\\
     & s_{M-1} & s_1  & s_2 & \ddots&  \vspace{0.2 cm}\\
    \vdots & \ddots & \ddots & \ddots & \ddots &s_3 \vspace{0.4 cm}\\
    && &  & &s_2\vspace{0.2 cm}\\
    s_2& \ldots&  && s_{M-1}  &  s_1\vspace{0.2 cm}\\
  \end{array}
\right),
$$
where each row is a cyclic shift of the preceding row, then matrix
$S$ has eigenvector
\begin{equation*}
\displaystyle
 y^{(j)}= \frac{1}{\sqrt{M-1}}\left(
\exp\left(-\frac{2\pi ji}{M-1}\right),\;\cdots,\; \exp\left(-\frac{2\pi
j(M-2)i}{M-1}\right),1\right)^{T},
\end{equation*}
and the corresponding eigenvalue
\begin{equation*}
\displaystyle \lambda_j(S)=\sum\limits_{\ell=1}^{M-1}s_{\ell}
\exp\left(-\frac{2\pi j\ell
i}{M-1}\right),\;\;i=\sqrt{-1},\;\;j=1,\ldots,M-1.
\end{equation*}}

{{\it{\bf Theorem 1.} The difference equations (15) and (16) are both
 uniquely solvable.}}

\begin{proof} From Lemma 3, we know that the eigenvalues of the
matrices
$\left(A-g_{0}^{(\alpha,\beta)}B\right)$
and
$\left(\widetilde{A}-g_{0}^{(\alpha,\beta)}\widetilde{B}\right)$
are
$$
\begin{array}{lll}
\displaystyle\lambda_j=\displaystyle\left[1-\frac{8}{45}\sin^4\left(
\frac{\pi
j}{M-1}\right)\right]+4g_{0}^{(\alpha,\beta)}
\sin^2\left( \frac{\pi j}{M-1}\right)\left[1+\frac{1}{3}\sin^2\left( \frac{\pi
j}{M-1}\right)\right],\;\;\;j=1,\cdots,M-1,
\end{array}
$$
and
$$
\begin{array}{lll}
\displaystyle\widetilde{\lambda}_j=&\displaystyle\left[1-\frac{4}{35}\sin^6\left(
\frac{\pi
j}{M-1}\right)\right]+4g_{0}^{(\alpha,\beta)}
\sin^2\left( \frac{\pi j}{M-1}\right)\vspace{0.3
cm}\\&\displaystyle\times\left[1+\frac{1}{3}\sin^2\left( \frac{\pi
j}{M-1}\right)+\frac{8}{45}\sin^4\left( \frac{\pi
j}{M-1}\right)\right],\;\;\;j=1,\cdots,M-1,
\end{array}
$$
respectively.

 Note that $\displaystyle\mu_{\alpha},\;\mu_{\beta}>0$ and
$g_{0}^{(1-\alpha)},\;g_{0}^{(1-\beta)}>0$, $\displaystyle
{\lambda}_j,\;\widetilde{\lambda}_j>0$. Thus
\begin{equation*}
\displaystyle \det\left(A-g_{0}^{(\alpha,\beta)}B\right)=\mathop{\prod}_{j=1}^{M-1}{\lambda}_j>0
\end{equation*}
and
\begin{equation*}
\displaystyle
\det\left(\widetilde{A}-g_{0}^{(\alpha,\beta)}\widetilde{B}\right)=\mathop{\prod}_{j=1}^{M-1}{\widetilde{\lambda}}_j>0.
\end{equation*}
Therefore, the above two matrices are both nonsingular. The
difference equations (13) and (14) are uniquely solvable. The proof is complete.
\end{proof}

\section{Stability Analysis }
 \quad\; In this section, we analyze
the stability of the difference schemes (13) and (14) by using the
Fourier method.

\subsection{Stability Analysis of Numerical Scheme (13)}

{\it{\bf Lemma 4 (\cite{CLTA,LD}).} The coefficients
$\displaystyle\varpi_{\ell}^{(1-\gamma)}$
$\displaystyle(\ell=0,1,\cdots)$ satisfy
$$\displaystyle
\begin{array}{lll}
\displaystyle (i)\;\;\displaystyle
\varpi_{0}^{(1-\gamma)}=1,\;\;\varpi_{1}^{(1-\gamma)}=\gamma-1,\;\;\varpi_{\ell}^{(1-\gamma)}<0,\;\;
 \ell\geq1;\vspace{0.2 cm}\\
(ii)\;\;\displaystyle
\sum\limits_{\ell=0}^{\infty}\varpi_{\ell}^{(1-\gamma)}=0;\;\forall\;
k\in
\mathbb{N}^{+},-\sum\limits_{\ell=1}^{k}\varpi_{\ell}^{(1-\gamma)}<1.
\end{array}
$$}
{\it{\bf Lemma 5.} The coefficients $g_{\ell}^{(1-\gamma)}$
$\displaystyle(\ell=0,1,\cdots)$ satisfy
$$\displaystyle
\begin{array}{lll}
\displaystyle (i)\;\;\displaystyle
g_{0}^{(1-\gamma)}=\frac{2-\gamma}{2},\;\;g_{1}^{(1-\gamma)}
=\frac{-\gamma^2+4\gamma-2}{2},\;\;g_{\ell}^{(1-\gamma)}<0,\;\;
 \ell\geq2;\vspace{0.2 cm}\\
(ii)\;\;\displaystyle
\sum\limits_{\ell=0}^{\infty}g_{\ell}^{(1-\gamma)}=0;\;\forall\;
k\in
\mathbb{N}^{+},-\sum\limits_{\ell=1}^{k}g_{\ell}^{(1-\gamma)}<\frac{2-\gamma}{2}.
\end{array}
$$}

\begin{proof} (i) From the above analysis, we easily obtain the expressions of
$g_0^{(1-\gamma)}$, $g_1^{(1-\gamma)}$, and
$$\begin{array}{lll}
\displaystyle
g_\ell^{(1-\gamma)}&=&\displaystyle\frac{2-\gamma}{2}\varpi_\ell^{(1-\gamma)}
+\frac{\gamma}{2}\varpi_{\ell-1}^{(1-\gamma)}\vspace{0.2 cm}\\
&=&\displaystyle\frac{2-\gamma}{2}\varpi_\ell^{(1-\gamma)}
+\frac{\gamma\ell}{2(\ell+\gamma-2)}\varpi_{\ell}^{(1-\gamma)}\vspace{0.2 cm}\\
&=&\displaystyle
\frac{2\ell-(2-\gamma)^2}{2(\ell+\gamma-2)}\varpi_{\ell}^{(1-\gamma)}.
\end{array}
$$
One has $g_{\ell}^{(1-\gamma)}\leq0$ for $\ell\geq2$ if $0<\gamma<1$.

(ii) In view of Lemma 4, it is not difficult to obtain these
relations by direct computations.
\end{proof}

 Let $U_j^k$ be the
approximate solution of (13) and define
$$
\begin{array}{lll}
\displaystyle
\rho_j^k=u_j^k-U_j^k,\;\;j=1,2,\cdots,M-1,\;\;k=0,1,\cdots,N,
\end{array}
$$
and
$$
\begin{array}{lll}
\displaystyle
\rho^{k}=\left(\rho_1^k,\rho_2^k,\cdots,\rho_{M-1}^k\right)^T,\;\;k=0,1,\cdots,N,
\end{array}
$$
respectively.

So, we can easily get the following roundoff error
equation
$$
\begin{array}{lll}
\displaystyle
\left[-\frac{1}{90}+\frac{1}{12}g_{0}^{(\alpha,\beta)}\right]\rho_{j-2}^{k+1}
+\left[\frac{2}{45}-\frac{4}{3}g_{0}^{(\alpha,\beta)}\right]\rho_{j-1}^{k+1}
+\left[\frac{14}{15}+\frac{5}{2}g_{0}^{(\alpha,\beta)}\right]\rho_{j}^{k+1}\vspace{0.2 cm}\\ \displaystyle
+\left[\frac{2}{45}-\frac{4}{3}g_{0}^{(\alpha,\beta)}\right]\rho_{j+1}^{k+1}
\displaystyle+\left[-\frac{1}{90}+\frac{1}{12}g_{0}^{(\alpha,\beta)}\right]\rho_{j+2}^{k+1}=
\left[-\frac{1}{90}-\frac{1}{12}g_{1}^{(\alpha,\beta)}\right]\rho_{j-2}^{k}\vspace{0.2
cm}\\ \displaystyle\
+\left[\frac{2}{45}+\frac{4}{3}g_{1}^{(\alpha,\beta)}\right]\rho_{j-1}^{k}
+\left[\frac{14}{15}-\frac{5}{2}g_{1}^{(\alpha,\beta)}\right]\rho_{j}^{k}
+\left[\frac{2}{45}+\frac{4}{3}g_{1}^{(\alpha,\beta)}\right]\rho_{j+1}^{k}\vspace{0.2 cm}\\
\displaystyle
+\left[-\frac{1}{90}-\frac{1}{12}g_{1}^{(\alpha,\beta)}\right]\rho_{j+2}^{k}
\displaystyle-\frac{1}{12} \sum\limits_{\ell=2}^{k+1}
g_{\ell}^{(\alpha,\beta)}\rho_{j-2}^{k+1-\ell}+\frac{4}{3}
\sum\limits_{\ell=2}^{k+1} g_{\ell}^{(\alpha,\beta)}\rho_{j-1}^{k+1-\ell}\vspace{0.2
cm}\\ \displaystyle-\frac{5}{2} \sum\limits_{\ell=2}^{k+1}
g_{\ell}^{(\alpha,\beta)}\rho_{j}^{k+1-\ell}
+\frac{4}{3} \sum\limits_{\ell=2}^{k+1}
g_{\ell}^{(\alpha,\beta)}\rho_{j+1}^{k+1-\ell}-\frac{1}{12} \sum\limits_{\ell=2}^{k+1}
g_{\ell}^{(\alpha,\beta)}\rho_{j+2}^{k+1-\ell}
,\vspace{0.2 cm}\\
\displaystyle\hspace{4.5 cm}j=1,2,\cdots,M-1,\;k=0,1,\cdots,\,N-1.
\end{array}\eqno(22)
$$
$$
\begin{array}{lll}
\displaystyle \rho_{0}^{k}=\rho_{M}^{k}=0,\;k=0,1,\cdots,\,N.
\end{array}
$$

Now, we define the grid functions
$$\displaystyle \rho^{k}(x)=\left\{
\begin{array}{lll}
\displaystyle{\rho_{j}^{k},\;\; \textrm{when}\;\;\; x_j-\frac{h}{2}<
x \leq x_j+\frac{h}{2},\;\;\;j=1,2,\cdots,M-1,}\vspace{0.3 cm}\\0,\;\;\;\;
\displaystyle\textrm{when}\;\; -h\leq x \leq
\frac{h}{2}\;\;\textrm{or}\;\; L-\frac{h}{2}<x\leq L+h,
\end{array}\right.
$$
then $\displaystyle \rho^k\left(x\right)$ can be expanded in a
Fourier series
$$
\begin{array}{lll}
\displaystyle
\rho^{k}\left(x\right)=\sum\limits_{l=-\infty}^{\infty}\xi_{k}\left(l\right)\exp\left(\frac{2\pi
l x}{L}i\right),
\end{array}
$$
where
$$
\begin{array}{lll}
\displaystyle
\xi_k\left(l\right)=\frac{1}{L}\int_{0}^{L}\rho^{k}\left(x\right)\exp\left(-\frac{2\pi
l x}{L}i\right)dx.
\end{array}
$$
Let
$$
\begin{array}{lll}
\displaystyle
\left\|\rho^k\right\|_2=\left(\sum\limits_{j=1}^{M-1}h\left|\rho_{j}^{k}\right|^{2}\right)^{\frac{1}{2}}=
\left[\int_{0}^{L}\left|\rho^k(x)\right|^2dx\right]^{\frac{1}{2}}.
\end{array}
$$
By the Parseval equality
$$
\begin{array}{lll}
\displaystyle
\int_{0}^{L}\left|\rho^k(x)\right|^2dx=\sum\limits_{l=-\infty}^{\infty}\left|\xi_k(l)\right|^2,
\end{array}
$$
one has
$$
\begin{array}{lll}
\displaystyle
\left\|\rho^k\right\|_2^2=\sum\limits_{l=-\infty}^{\infty}\left|\xi_k(l)\right|^2.
\end{array}
$$

Now we suppose that the solution of
equation (22) has the following form
$$
\begin{array}{lll}
\displaystyle \rho_{j}^{k}=\xi_{k} \exp\left(i\beta j h \right),
\end{array}
$$
where $\displaystyle \beta={2\pi l}/{L}$. \\
Substituting the above expression into (22) gives
$$
\begin{array}{lll}
\displaystyle \mathcal {Q}\xi_{k+1}=\displaystyle \mathcal
{P}\xi_{k} -4\sin^2\left( \frac{\beta
h}{2}\right)\left[1+\frac{1}{3}\sin^2\left( \frac{\beta
h}{2}\right)\right]
\sum\limits_{\ell=2}^{k+1}
g_{\ell}^{(\alpha,\beta)}
\xi_{k+1-\ell},\;\;k=0,1,\ldots,N-1.
\end{array}\eqno(23)
$$
where
$$
\begin{array}{lll}
\displaystyle\mathcal {Q}=\left[1-\frac{8}{45}\sin^4\left(
\frac{\beta
h}{2}\right)\right]+4g_{0}^{(\alpha,\beta)}
\sin^2\left( \frac{\beta h}{2}\right)\left[1+\frac{1}{3}\sin^2\left(
\frac{\beta h}{2}\right)\right],
\end{array}
$$
$$
\begin{array}{lll}
\displaystyle\mathcal {P}=\left[1-\frac{8}{45}\sin^4\left(
\frac{\beta
h}{2}\right)\right]-4g_{1}^{(\alpha,\beta)}
\sin^2\left( \frac{\beta h}{2}\right)\left[1+\frac{1}{3}\sin^2\left(
\frac{\beta h}{2}\right)\right].
\end{array}
$$

 {\it{\bf Lemma 5.} If
$\mathcal {Q}$ and $\mathcal {P}$ are defined as above,
 then
$$
\begin{array}{lll}
\displaystyle \left| \frac{\mathcal {P}}{\mathcal {Q}}\right|\leq1,
\end{array}
$$}
\begin{proof}One can show that
$$
\begin{array}{lll}
\displaystyle (\mathcal {P}+\mathcal {Q})(\mathcal {P}-\mathcal
{Q})&=&\displaystyle-16\left[\mu_{\alpha}\left(g_{0}^{(1-\alpha)}-g_{1}^{(1-\alpha)}\right)
+\mu_{\beta}\left(g_{0}^{(1-\beta)}-g_{1}^{(1-\beta)}\right)\right]
\vspace{0.2cm}\\&&\displaystyle\times
\left[\mu_{\alpha}\left(g_{0}^{(1-\alpha)}+g_{1}^{(1-\alpha)}\right)
+\mu_{\beta}\left(g_{0}^{(1-\beta)}+g_{1}^{(1-\beta)}\right)\right]
\vspace{0.2cm}\\&&\displaystyle\times \sin^4\left( \frac{\beta
h}{2}\right)\left[1+\frac{1}{3}\sin^2\left( \frac{\beta
h}{2}\right)\right]^2\vspace{0.2cm}\\
&=&\displaystyle-16\left[\frac{(1-\alpha)(4-\alpha)}{2}\mu_{\alpha}
+\frac{(1-\beta)(4-\beta)}{2}\mu_{\beta}\right]
\vspace{0.2cm}\\&&\displaystyle\times
\left[\frac{\alpha(3-\alpha)}{2}\mu_{\alpha}
+\frac{\beta(3-\beta)}{2}\mu_{\beta}\right]\sin^4\left( \frac{\beta
h}{2}\right)\left[1+\frac{1}{3}\sin^2\left( \frac{\beta
h}{2}\right)\right]^2
\end{array}
$$
Note that $\displaystyle \mu_{\alpha},\;\mu_{\beta}>0,$ and
$\displaystyle 0<\alpha,\;\beta<1$, therefore we obtain that $\displaystyle
(\mathcal {P}+\mathcal {Q})(\mathcal {P}-\mathcal {Q})\leq0$, i.e.,
$$
\begin{array}{lll}
\displaystyle \left| \frac{\mathcal {P}}{\mathcal {Q}}\right|\leq 1.
\end{array}
$$
This ends the proof.
\end{proof}

 {\it{\bf Lemma 6.} If time and space steps $\tau$ and $h$
 satisfy
$$
\begin{array}{lll}
\displaystyle
\frac{\tau^{\alpha}\left(-\alpha^2+4\alpha-2\right)\mathcal
{A}+\tau^{\beta}\left(-\beta^2+4\beta-2\right)\mathcal
{B}}{h^2}\leq\frac{37}{120},
\end{array}\eqno(24)
$$
then one has
$$
\begin{array}{lll}
\displaystyle \mathcal {P}\geq0.
\end{array}
$$}

\begin{proof} If $\tau$ and $h$ satisfy
$$
\begin{array}{lll}
\displaystyle
\frac{\tau^{\alpha}\left(-\alpha^2+4\alpha-2\right)\mathcal
{A}+\tau^{\beta}\left(-\beta^2+4\beta-2\right)\mathcal
{B}}{h^2}\leq0,
\end{array}
$$
 we easily obtain $\mathcal {P}\geq0$.

 In effect,
$$
\begin{array}{lll}
&\displaystyle
0\leq\frac{\tau^{\alpha}\left(-\alpha^2+4\alpha-2\right)\mathcal
{A}+\tau^{\beta}\left(-\beta^2+4\beta-2\right)\mathcal
{B}}{h^2}\leq\frac{37}{120}\vspace{0.2 cm}\\
\Rightarrow&\displaystyle
0\leq\frac{16}{3}g_{1}^{(\alpha,\beta)}\leq
1-\frac{8}{45}\sin^4\left( \frac{\beta h}{2}\right)\vspace{0.2 cm}\\
\Rightarrow&\displaystyle
0\leq4g_{1}^{(\alpha,\beta)}
\sin^2\left( \frac{\beta h}{2}\right)\left[1+\frac{1}{3}\sin^2\left(
\frac{\beta h}{2}\right)\right]\leq 1-\frac{8}{45}\sin^4\left(
\frac{\beta h}{2}\right).
\end{array}
$$
It immediately follows that
$$
\begin{array}{lll}
\displaystyle\mathcal {P}\geq0.
\end{array}
$$
The proof is complete.
\end{proof}

{\it{\bf Lemma 7.} Suppose that $\xi_{k+1}$ $(k=0,1,\cdots,N-1)$
is the solution of equation (23). Under the condition of (24), it follows that
$$
\begin{array}{lll}\displaystyle
\left|\xi_{k+1}\right|\leq \left|\xi_0\right|,
\;\;\;k=0,1,\cdots,N-1.
\end{array}
$$}
\begin{proof} For $k=0$, from equation (23), we have
$$
\begin{array}{lll}
\displaystyle |\xi_1|=\left|\frac{\mathcal {P}}{ \mathcal
{Q}}\right||\xi_0|.
\end{array}
$$

According to Lemma 5 it is clear
 that
$$
\begin{array}{lll}
\displaystyle |\xi_1|\leq |\xi_0|.
\end{array}
$$

Now, we suppose that
$$
\begin{array}{lll}
\displaystyle |\xi_\ell|\leq|\xi_0|,\;\;\;(\ell=1,2,\cdots,k).
\end{array}
$$

For $k>0$, from equation (23), Lemmas 4 and 5,and the condition
of Lemma 6, i.e., $\mathcal {P}\geq0$, we have
$$
\begin{array}{lll}
\displaystyle \mathcal {Q}|\xi_{k+1}|&=&\displaystyle \left|\mathcal
{P}\xi_{k} -4\sin^2\left( \frac{\beta
h}{2}\right)\left[1+\frac{1}{3}\sin^2\left( \frac{\beta
h}{2}\right)\right]\sum\limits_{\ell=2}^{k+1}
g_{\ell}^{(\alpha,\beta)}
\xi_{k+1-\ell}\right|\vspace{0.2 cm}\\&\leq& \displaystyle
\left|\mathcal {P}\right|\left|\xi_{k}\right|+4\sin^2\left(
\frac{\beta h}{2}\right)\left[1+\frac{1}{3}\sin^2\left( \frac{\beta
h}{2}\right)\right]
\sum\limits_{\ell=2}^{k+1}\left|
g_{\ell}^{(\alpha,\beta)}\right|
\left|\xi_{k+1-\ell}\right|\vspace{0.2 cm}\\&\leq& \displaystyle
\left\{\left|\mathcal {P}\right|+4\sin^2\left( \frac{\beta
h}{2}\right)\left[1+\frac{1}{3}\sin^2\left( \frac{\beta
h}{2}\right)\right] \sum\limits_{\ell=2}^{k+1}\left|
g_{\ell}^{(\alpha,\beta)}\right|
\right\}\left|\xi_{0}\right|\vspace{0.2 cm}\\&\leq& \displaystyle
\left\{\mathcal {P}-4\sin^2\left( \frac{\beta
h}{2}\right)\left[1+\frac{1}{3}\sin^2\left( \frac{\beta
h}{2}\right)\right] \left[
g_{0}^{(\alpha,\beta)}+
g_{1}^{(\alpha,\beta)}
\right] \right\}\left|\xi_{0}\right|\vspace{0.2 cm}\\&=&\mathcal {Q}
\left|\xi_{0}\right|,
\end{array}
$$
that is,
$$
\begin{array}{lll}
\displaystyle\left|\xi_{k+1}\right|\leq\left|\xi_0\right|.
\end{array}
$$
The proof is thus completed.
\end{proof}

 {\it{\bf Theorem 3.}
Under condition (24), the difference scheme (13) is stable.}

\begin{proof} According to Lemma 7, we obtain
$$
\begin{array}{lll}
\displaystyle
\|\rho^{k+1}\|_2&=&\displaystyle\left(\sum\limits_{j=1}^{M-1}h\left|\rho_{j}^{k+1}\right|^{2}\right)^{\frac{1}{2}}
=\left(\sum\limits_{j=1}^{M-1}h\left|\xi_{k+1} \exp\left(i\beta j h
\right)\right|^{2}\right)^{\frac{1}{2}}
=\left(\sum\limits_{j=1}^{M-1}h\left|\xi_{k+1}\right|^{2}\right)^{\frac{1}{2}}
\vspace{0.2 cm}\\
&\leq&\displaystyle\left(\sum\limits_{j=1}^{M-1}h\left|\xi_{0}\right|^{2}\right)^{\frac{1}{2}}
=\left(\sum\limits_{j=1}^{M-1}h\left|\xi_{0} \exp\left(i\beta j h
\right)\right|^{2}\right)^{\frac{1}{2}}
=\left(\sum\limits_{j=1}^{M-1}h\left|\rho_{j}^{0}\right|^{2}\right)^{\frac{1}{2}}
\vspace{0.3 cm}\\
&=&\displaystyle\|\rho^0\|_2,\;\;k=0,1,\cdots,N-1,
\end{array}
$$
 which means that the difference
scheme (13) is stable. The proof is complete.
\end{proof}

\subsection{Stability Analysis of Numerical Scheme (14)}

Similarly,  let $\widetilde{U}_j^k$ be the approximate solution of
(14) and define $$
\begin{array}{lll}
\displaystyle
\widetilde{\rho}_j^k=u_j^k-\widetilde{U}_j^k,\;\;j=1,2,\cdots,M-1,\;\;k=0,1,\cdots,N,
\end{array}
$$
then we can get truncation error equation of (14) which is
$$
\begin{array}{lll}
\displaystyle
\left[\frac{1}{560}-\frac{1}{90}g_{0}^{(\alpha,\beta)}\right]\widetilde{\rho}_{j-3}^{k+1}+
\left[-\frac{3}{280}+\frac{3}{20}g_{0}^{(\alpha,\beta)}\right]\widetilde{\rho}_{j-2}^{k+1}
+\left[\frac{3}{112}-\frac{3}{2}g_{0}^{(\alpha,\beta)}\right]\widetilde{\rho}_{j-1}^{k+1}\vspace{0.2
cm}\\ \displaystyle\
+\left[\frac{27}{28}+\frac{49}{18}g_{0}^{(\alpha,\beta)}\right]\widetilde{\rho}_{j}^{k+1}
+\left[\frac{3}{112}-\frac{3}{2}g_{0}^{(\alpha,\beta)}\right]\widetilde{\rho}_{j+1}^{k+1}
\displaystyle+\left[-\frac{3}{280}+\frac{3}{20}g_{0}^{(\alpha,\beta)}\right]\widetilde{\rho}_{j+2}^{k+1}
\vspace{0.2 cm}\\\displaystyle+
\left[\frac{1}{560}-\frac{1}{90}g_{0}^{(\alpha,\beta)}\right]\widetilde{\rho}_{j+3}^{k+1}
=\left[\frac{1}{560}+\frac{1}{90}g_{1}^{(\alpha,\beta)}\right]\widetilde{\rho}_{j-3}^{k}
-\left[\frac{3}{280}+\frac{3}{20}g_{1}^{(\alpha,\beta)}\right]\widetilde{\rho}_{j-2}^{k}\vspace{0.2
cm}\\ \displaystyle\
+\left[\frac{3}{112}+\frac{3}{2}g_{1}^{(\alpha,\beta)}\right]\widetilde{\rho}_{j-1}^{k}
+\left[\frac{27}{28}-\frac{49}{18}g_{1}^{(\alpha,\beta)}\right]\widetilde{\rho}_{j}^{k}
+\left[\frac{3}{112}+\frac{3}{2}g_{1}^{(\alpha,\beta)}\right]\widetilde{\rho}_{j+1}^{k}
\vspace{0.2 cm}\\\displaystyle
-\left[\frac{3}{280}+\frac{3}{20}g_{1}^{(\alpha,\beta)}\right]\widetilde{\rho}_{j+2}^{k}
+\left[\frac{1}{560}+\frac{1}{90}g_{1}^{(\alpha,\beta)}\right]\widetilde{\rho}_{j+3}^{k}
 +\frac{1}{90}
\sum\limits_{\ell=2}^{k+1}g_{\ell}^{(\alpha,\beta)}\widetilde{\rho}_{j-3}^{k+1-\ell}\vspace{0.2 cm}\\\displaystyle
-\frac{3}{20} \sum\limits_{\ell=2}^{k+1}
 g_{\ell}^{(\alpha,\beta)}\widetilde{\rho}_{j-2}^{k+1-\ell}
 +\frac{3}{2}
\sum\limits_{\ell=2}^{k+1}g_{\ell}^{(\alpha,\beta)}\widetilde{\rho}_{j-1}^{k+1-\ell}
-\frac{49}{18} \sum\limits_{\ell=2}^{k+1}
 g_{\ell}^{(\alpha,\beta)}\widetilde{\rho}_{j}^{k+1-\ell}
\vspace{0.2 cm}\\
\displaystyle +\frac{3}{2} \sum\limits_{\ell=2}^{k+1}
 g_{\ell}^{(\alpha,\beta)}\widetilde{\rho}_{j+1}^{k+1-\ell}
-\frac{3}{20} \sum\limits_{\ell=2}^{k+1}
 g_{\ell}^{(\alpha,\beta)}\widetilde{\rho}_{j+2}^{k+1-\ell}
 +\frac{1}{90}
\sum\limits_{\ell=2}^{k+1}g_{\ell}^{(\alpha,\beta)}\widetilde{\rho}_{j+3}^{k+1-\ell}
,\vspace{0.2 cm}\\\displaystyle\hspace{4.5 cm} j=1,2,\cdots,M-1,\;k=0,1,\cdots,N-1.
\end{array}\eqno(25)
$$
$$
\begin{array}{lll}
\displaystyle
\widetilde{\rho}_{0}^{k}=\widetilde{\rho}_{M}^{k}=0,\;k=0,1,\cdots,N.
\end{array}
$$

Define the grid functions as
$$\displaystyle \widetilde{\rho}^{k}(x)=\left\{
\begin{array}{lll}
\displaystyle{\widetilde{\rho}_{j}^{k},\;\; \textrm{when}\;\;\;
x_j-\frac{h}{2}< x \leq
x_j+\frac{h}{2},\;\;j=1,2,\cdots,M-1,}\vspace{0.3 cm}\\0,\;\;\;\;
\displaystyle\textrm{when}\;\; -2h\leq x \leq
\frac{h}{2}\;\;\textrm{or}\;\; L-\frac{h}{2}<x\leq L+2h.
\end{array}\right.
$$
The function $\displaystyle \widetilde{\rho}^k\left(x\right)$ can be
expanded in a Fourier series
$$
\begin{array}{lll}
\displaystyle
\widetilde{\rho}^{k}\left(x\right)=\sum\limits_{l=-\infty}^{\infty}\widetilde{\xi}_{k}\left(l\right)\exp\left(\frac{2\pi
l x}{L}i\right),
\end{array}
$$
where
$$
\begin{array}{lll}
\displaystyle
\widetilde{\xi}_k\left(l\right)=\frac{1}{L}\int_{0}^{L}\widetilde{\rho}^{k}\left(x\right)\exp\left(-\frac{2\pi
l x}{L}i\right)dx,\;i^2=-1.
\end{array}
$$

Letting
$$
\begin{array}{lll}
\displaystyle \widetilde{\rho}_{j}^{k}=\widetilde{\xi}_{k}
\exp\left(i\beta j h \right),
\end{array}
$$
and substituting it into (25) yield
$$
\begin{array}{lll}
\displaystyle \mathcal {\widetilde{Q}}\xi_{k+1}=&\displaystyle
\mathcal {\widetilde{P}}\widetilde{\xi}_{k} -4\sin^2\left(
\frac{\beta h}{2}\right)\left[1+\frac{1}{3}\sin^2\left( \frac{\beta
h}{2}\right)+\frac{8}{45}\sin^4\left( \frac{\beta
h}{2}\right)\right]\vspace{0.2 cm}\\& \displaystyle\times
\sum\limits_{\ell=2}^{k+1}
g_{\ell}^{(\alpha,\beta)}
\widetilde{\xi}_{k+1-\ell},\;\;k=0,1,\ldots,N-1,
\end{array}\eqno(26)
$$
where
$$
\begin{array}{lll}
\displaystyle\mathcal
{\widetilde{Q}}=&\displaystyle\left[1-\frac{4}{35}\sin^6\left(
\frac{\beta
h}{2}\right)\right]+4g_{0}^{(\alpha,\beta)}
\sin^2\left( \frac{\beta h}{2}\right)\left[1+\frac{1}{3}\sin^2\left( \frac{\beta
h}{2}\right)+\frac{8}{45}\sin^4\left( \frac{\beta
h}{2}\right)\right],
\end{array}
$$
$$
\begin{array}{lll}
\displaystyle\mathcal
{\widetilde{P}}=&\displaystyle\left[1-\frac{4}{35}\sin^6\left(
\frac{\beta
h}{2}\right)\right]-4g_{1}^{(\alpha,\beta)}
\sin^2\left( \frac{\beta h}{2}\right)\left[1+\frac{1}{3}\sin^2\left( \frac{\beta
h}{2}\right)+\frac{8}{45}\sin^4\left( \frac{\beta
h}{2}\right)\right].
\end{array}
$$

The following lemmas and theorem can be similarly proved.

 {\it{\bf Lemma 8.} If
$\mathcal {\widetilde{Q}}$ and $\mathcal {\widetilde{P}}$ are
defined as above,
 then
$$
\begin{array}{lll}
\displaystyle \left| \frac{\mathcal {\widetilde{P}}}{\mathcal
{\widetilde{Q}}}\right|\leq1,
\end{array}
$$}

{\it{\bf Lemma 9.} If time and space steps $\tau$ and $h$
 satisfy
$$
\begin{array}{lll}
\displaystyle
\frac{\tau^{\alpha}\left(-\alpha^2+4\alpha-2\right)\mathcal
{A}+\tau^{\beta}\left(-\beta^2+4\beta-2\right)\mathcal
{B}}{h^2}\leq\frac{279}{952},
\end{array}\eqno(27)
$$
then
$$
\begin{array}{lll}
\displaystyle \mathcal {\widetilde{P}}\geq0.
\end{array}
$$}

{\it{\bf Lemma 10.} Supposing that $\widetilde{\xi}_{k+1}$
$(k=0,1,\cdots,N-1)$ \red{is} the solution of equation (26), under
condition (27), then it follows that
$$
\begin{array}{lll}\displaystyle
\left|\widetilde{\xi}_{k+1}\right|\leq
\left|\widetilde{\xi}_0\right|, \;\;\;k=0,1,\cdots,N-1.
\end{array}
$$}

 {\it{\bf Theorem 4.}
Under condition (27), the difference scheme (14) is stable.}

\section { Convergence Analysis }

\qquad In this section, we study the convergence of schemes (13) and (14).

\subsection { Convergence Analysis of Numerical Scheme (13)}

For equation (13), suppose that
$${E}_{j}^{k}=u(x_j,t_k)-u_j^{k},\;\;j=1,\cdots,M-1,k=1,\cdots,N,$$
and denote
$$E^{k}=\left(E_{1}^{k},
E_{2}^{k},\cdots,E_{M-1}^{k}\right)^{T},\;\;
{R}^{k}=\left(R_{1}^{k},
R_{2}^{k},\cdots,R_{M-1}^{k}\right)^{T},\;\;k=1,\cdots,N.
$$

Then we obtain
$$
\begin{array}{lll} \displaystyle
 \frac{E_j^{k+1}-E_j^k}{\tau}
 = \frac{\mathcal {A}}{\tau^{1-\alpha}}\mathscr{L}_1
\sum\limits_{\ell=0}^{k+1}g_\ell^{(1-\alpha)}E_j^{k+1-\ell}
+\frac{\mathcal {B}}{\tau^{1-\beta}}\mathscr{L}_1
\sum\limits_{\ell=0}^{k+1}g_\ell^{(1-\beta)}E_j^{k+1-\ell}
\vspace{0.2 cm}\\
\displaystyle
  +f_j^{k+\frac{1}{2}}+R_j^{k+1},\;\;\;
0\leq k\leq N-1,\;1\leq j\leq M-1.
\end{array}\eqno(28)
$$
Similar to the stability analysis above, we define the grid functions
$$\displaystyle E^{k}(x)=\left\{
\begin{array}{lll}
\displaystyle{E_{j}^{k},\;\; \textrm{when}\;\;\; x_j-\frac{h}{2}< x
\leq x_j+\frac{h}{2},\;\;j=1,2,\cdots,M-1,}\vspace{0.3 cm}\\0,\;
\displaystyle\textrm{when}\;\;\; -h\leq x \leq \frac{h}{2}
\;\;\textrm{or}\;\; L-\frac{h}{2}<x\leq L+h,
\end{array}\right.
$$
and
$$\displaystyle {R}^{k}(x)=\left\{
\begin{array}{lll}
\displaystyle{R_{j}^{k},\;\; \textrm{when}\;\;\; x_j-\frac{h}{2}< x
\leq x_j+\frac{h}{2},\;\;j=1,2,\cdots,M-1,}\vspace{0.3 cm}\\0,\;
\displaystyle\textrm{when}\;\;\; -h\leq x \leq \frac{h}{2}
\;\;\textrm{or}\;\; L-\frac{h}{2}<x\leq L+h.
\end{array}\right.
$$
Functions $\displaystyle E^{k}(x)$ and $\displaystyle{R}^{k}(x)$ can be
expanded into the following Fourier series, respectively,
$$
\begin{array}{lll}
\displaystyle
E^{k}\left(x\right)=\sum\limits_{l=-\infty}^{\infty}\zeta_{k}\left(l\right)\exp\left(\frac{2\pi
l x}{L}i\right),
\end{array}
$$
and
$$
\begin{array}{lll}
\displaystyle
{R}^{k}\left(x\right)=\sum\limits_{l=-\infty}^{\infty}\eta_{k}\left(l\right)\exp\left(\frac{2\pi
l x}{L}i\right),
\end{array}
$$
where
$$
\begin{array}{lll}
\displaystyle
\zeta_k\left(l\right)=\frac{1}{L}\int_{0}^{L}E^{k}\left(x\right)\exp\left(-\frac{2\pi
l x}{L}i\right)dx,
\end{array}
$$
and
$$
\begin{array}{lll}
\displaystyle
\eta_k\left(l\right)=\frac{1}{L}\int_{0}^{L}{R}^{k}\left(x\right)\exp\left(-\frac{2\pi
l x}{L}i\right)dx.
\end{array}
$$

The 2-norms are given below
$$
\begin{array}{lll}
\displaystyle
\|E^{k}\|_2&=&\displaystyle\left(\sum\limits_{i=1}^{M-1}h\left|E_{i}^{k}\right|^{2}\right)^{\frac{1}{2}}
=\left(\sum\limits_{l=-\infty}^{\infty}\left|\zeta_k(l)\right|^2\right)^{\frac{1}{2}},
\end{array}\eqno(29)
$$
and
$$
\begin{array}{lll}
\displaystyle
\|{R}^{k}\|_2&=&\displaystyle\left(\sum\limits_{i=1}^{M-1}h\left|R_{i}^{k}\right|^{2}\right)^{\frac{1}{2}}
=\left(\sum\limits_{l=-\infty}^{\infty}\left|\eta_k(l)\right|^2\right)^{\frac{1}{2}}.
\end{array}\eqno(30)
$$

Assume that $E_{i}^{k}$ and
$R_{i}^{k}$ have the following forms
$$
\begin{array}{lll}
\displaystyle E_{j}^{k}=\zeta_{k} \exp\left(i\beta j h \right),
\end{array}
$$
and
$$
\begin{array}{lll}
\displaystyle R_{j}^{k}=\eta_{k} \exp\left(i\beta j h \right),
\end{array}
$$
respectively. Substituting the above two expressions into (28) yields
$$
\begin{array}{lll}
\displaystyle \mathcal {Q}\zeta_{k+1}=&\displaystyle \mathcal
{P}\zeta_{k} -4\sin^2\left( \frac{\beta
h}{2}\right)\left[1+\frac{1}{3}\sin^2\left( \frac{\beta
h}{2}\right)\right]\vspace{0.2 cm}\\& \displaystyle\times
\sum\limits_{\ell=2}^{k+1}
g_{\ell}^{(\alpha,\beta)}
\zeta_{k+1-\ell}+\tau\left[1-\frac{8}{45}\sin^4\left( \frac{\beta
h}{2}\right)\right]\eta_{k+1}.
\end{array}\eqno(31)
$$

{\it{\bf Lemma 11.} Let $\zeta_{k+1}$ $(k=0,1,\cdots,N-1)$ be the
solution of equation (31), under condition (24), then there
exists a positive constant $C_2$ such that
$$
\begin{array}{lll}\displaystyle
\left|\zeta_{k+1}\right|\leq C_{2}(k+1)\tau\left|\eta_1\right|,
\;\;\;k=0,1,\cdots,N-1.
\end{array}
$$}

\begin{proof} From $E^{0}=0$, we have
$$\zeta_0=\zeta_0(l)=0.$$

In addition, we know that there
exists a positive constant $C_1$ such that
$$
\begin{array}{lll}
\displaystyle \left|R_{j}^{k+1}\right| \leq C_{1}(\tau^2+h^6),
\end{array}
$$
and
$$
\begin{array}{lll}
\displaystyle \|R_{j}^{k+1}\|\leq C_{1}
\sqrt{(M-1)h}\left(\tau^2+h^6\right)\leq
C_{1}\sqrt{L}\left(\tau^2+h^6\right).
\end{array}
$$

In view of the convergence of series (30), there exists a
positive constant $C_2$ such that
$$
\begin{array}{lll}
\displaystyle \left|\eta_{k+1}\right|=\left|\eta_{k+1}(l)\right|
\leq C_{2}
\left|\eta_{1}\right|= C_{2}\left|\eta_{1}(l)\right|.\\
\end{array}\eqno(32)
$$

For $k=0$, from (31) we have
$$
\begin{array}{lll}
\displaystyle \zeta_{1}=\frac{ \tau}{\mathcal {Q}}\left[1-\frac{8}{45}\sin^4\left(
\frac{\beta h}{2}\right)\right]\eta_{k+1}.
\end{array}
$$
Note from equation (32) that one has
$$
\begin{array}{lll}
\displaystyle \left|\zeta_{1}\right| \leq \tau
\left|\eta_{1}\right|\leq C_{2}\tau \left|\eta_{1}\right|.
\end{array}
$$
Now, we suppose that
$$
\begin{array}{lll}
\displaystyle \left|\zeta_{\ell}\right| \leq C_{2}\ell\tau
\left|\eta_{1}\right|,\;\;\ell=1,\cdots,N-1.
\end{array}
$$
For $k>0$ and (24), one gets
$$
\begin{array}{lll}
\displaystyle \mathcal {Q}|\zeta_{k+1}|&=&\displaystyle
\left|\mathcal {P}\zeta_{k} -4\sin^2\left( \frac{\beta
h}{2}\right)\left[1+\frac{1}{3}\sin^2\left( \frac{\beta
h}{2}\right)\right]\right.\vspace{0.2 cm}\\&&
\displaystyle\left.\times \sum\limits_{\ell=2}^{k+1}
g_{\ell}^{(\alpha,\beta)}
\zeta_{k+1-\ell}+\tau\left[1-\frac{8}{45}\sin^4\left( \frac{\beta
h}{2}\right)\right]\eta_{k+1}\right|\vspace{0.2 cm}\\
&\leq&\displaystyle \mathcal {P}\left|\zeta_{k}\right|
+4\sin^2\left( \frac{\beta
h}{2}\right)\left[1+\frac{1}{3}\sin^2\left( \frac{\beta
h}{2}\right)\right]\vspace{0.2 cm}\\&& \displaystyle\times
\sum\limits_{\ell=2}^{k+1}
\left|g_{\ell}^{(\alpha,\beta)}\right|
\left|\zeta_{k+1-\ell}\right|+\tau\left[1-\frac{8}{45}\sin^4\left(
\frac{\beta h}{2}\right)\right]\left|\eta_{k+1}\right|\vspace{0.2 cm}\\
&\leq&\displaystyle \left\{\mathcal {P}k -4\sin^2\left( \frac{\beta
h}{2}\right)\left[1+\frac{1}{3}\sin^2\left( \frac{\beta
h}{2}\right)\right]\right.\vspace{0.2 cm}\\&&
\displaystyle\left.\times \sum\limits_{\ell=2}^{k+1}
g_{\ell}^{(\alpha,\beta)}
(k+1-\ell)+\left[1-\frac{8}{45}\sin^4\left( \frac{\beta
h}{2}\right)\right]\right\}C_2\tau\left|\eta_{1}\right|\vspace{0.2 cm}\\
&\leq&\displaystyle \left\{\mathcal {P}k -4\sin^2\left( \frac{\beta
h}{2}\right)\left[1+\frac{1}{3}\sin^2\left( \frac{\beta
h}{2}\right)\right]\right.\vspace{0.2 cm}\\&&
\displaystyle\left.\times \sum\limits_{\ell=2}^{\infty}
g_{\ell}^{(\alpha,\beta)}
k+\left[1-\frac{8}{45}\sin^4\left( \frac{\beta
h}{2}\right)\right]\right\}C_2\tau\left|\eta_{1}\right|\vspace{0.2 cm}\\
&\leq&\displaystyle\mathcal {Q}C_{2}(k+1)\tau\left|\eta_{1}\right|.
\end{array}
$$
Hence,
$$
\begin{array}{lll}
\displaystyle \left|\zeta_{k+1}\right| \leq
C_{2}(k+1)\tau\left|\eta_{1}\right|.
\end{array}
$$
The proof is completed.
\end{proof}

 {\it{\bf Theorem 5.} Under condition (24), the difference scheme (13) is convergent with
order $O(\tau^2+h^6)$.}

\begin{proof} Using (29), (30), Lemma 6, and
condition (24), one has
$$
\begin{array}{lll}
\displaystyle \left\|E^{k+1}\right\|_{2} \leq
C_{2}(k+1)\tau\left\|R^{1}\right\|_{2} \leq
C_{1}C_{2}\sqrt{L}(k+1)\tau\left(\tau^2+h^6\right) .
\end{array}
$$
Due to $k\leq N-1$, then
$$(k+1)\tau\leq T,$$
thus,
$$
\begin{array}{lll}
\displaystyle \left\|E^{k+1}\right\|_{2} \leq
C\left(\tau^2+h^6\right),
\end{array}
$$
where $C=C_1C_2T\sqrt{L}.$ This ends the proof.
\end{proof}

\subsection { Convergence Analysis of Numerical Scheme (14)}

Define
$${\widetilde{E}}_{i}^{k}=u(x_i,t_k)-u_i^{k},\;\;i=1,\cdots,M-1,k=1,\cdots,N,$$
and denote
$$\widetilde{E}^{k}=\left(\widetilde{E}_{1}^{k},
\widetilde{E}_{2}^{k},\cdots,\widetilde{E}_{M-1}^{k}\right)^{T},\;\;
{\widetilde{R}}^{k}=\left(\widetilde{R}_{1}^{k},
\widetilde{R}_{2}^{k},\cdots,\widetilde{R}_{M-1}^{k}\right)^{T},\;\;k=1,\cdots,N.
$$

From equation (14), one has
$$
\begin{array}{lll} \displaystyle
 \frac{\widetilde{E}_j^{k+1}-\widetilde{E}_j^k}{\tau}
 = \frac{\mathcal {A}}{\tau^{1-\alpha}}\mathscr{L}_2
\sum\limits_{\ell=0}^{k+1}g_\ell^{(1-\alpha)}\widetilde{E}_j^{k+1-\ell}
+\frac{\mathcal {B}}{\tau^{1-\beta}}\mathscr{L}_2
\sum\limits_{\ell=0}^{k+1}g_\ell^{(1-\beta)}\widetilde{E}_j^{k+1-\ell}
\vspace{0.2 cm}\\
\displaystyle
  +f_j^{k+\frac{1}{2}}+\widetilde{R}_j^{k+1},\;\;\;
0\leq k\leq N-1,\;1\leq j\leq M-1.
\end{array}\eqno(33)
$$
We now define the grid functions
$$\displaystyle \widetilde{E}^{k}(x)=\left\{
\begin{array}{lll}
\displaystyle{\widetilde{E}_{j}^{k},\;\; \textrm{when}\;\;\;
x_j-\frac{h}{2}< x \leq
x_j+\frac{h}{2},\;\;j=1,2,\cdots,M-1,}\vspace{0.3 cm}\\0,\;
\displaystyle\textrm{when}\;\;\; -2h\leq x \leq \frac{h}{2}
\;\;\textrm{or}\;\; L-\frac{h}{2}<x\leq L+2h,
\end{array}\right.
$$
and
$$\displaystyle {\widetilde{R}}^{k}(x)=\left\{
\begin{array}{lll}
\displaystyle{\widetilde{R}_{j}^{k},\;\; \textrm{when}\;\;\;
x_j-\frac{h}{2}< x \leq
x_j+\frac{h}{2},\;\;j=1,2,\cdots,M-1,}\vspace{0.3 cm}\\0,\;
\displaystyle\textrm{when}\;\;\; -2h\leq x \leq \frac{h}{2}
\;\;\textrm{or}\;\; L-\frac{h}{2}<x\leq L+2h.
\end{array}\right.
$$

The functions $\displaystyle \widetilde{E}^{k}(x)$ and
$\displaystyle{\widetilde{R}}^{k}(x)$ can be expanded into the following
Fourier series,
$$
\begin{array}{lll}
\displaystyle
\widetilde{E}^{k}\left(x\right)=\sum\limits_{l=-\infty}^{\infty}\widetilde{\zeta}_{k}\left(l\right)\exp\left(\frac{2\pi
l x}{L}i\right),
\end{array}
$$
and
$$
\begin{array}{lll}
\displaystyle
{\widetilde{R}}^{k}\left(x\right)=\sum\limits_{l=-\infty}^{\infty}\widetilde{\eta}_{k}\left(l\right)\exp\left(\frac{2\pi
l x}{L}i\right),
\end{array}
$$
where
$$
\begin{array}{lll}
\displaystyle
\widetilde{\zeta}_k\left(l\right)=\frac{1}{L}\int_{0}^{L}\widetilde{E}^{k}\left(x\right)\exp\left(-\frac{2\pi
l x}{L}i\right)dx,
\end{array}
$$
and
$$
\begin{array}{lll}
\displaystyle
\widetilde{\eta}_k\left(l\right)=\frac{1}{L}\int_{0}^{L}{\widetilde{R}}^{k}\left(x\right)\exp\left(-\frac{2\pi
l x}{L}i\right)dx.
\end{array}
$$

Similar to the above analysis, we assume that
$\widetilde{E}_{i}^{k}$ and $\widetilde{R}_{i}^{k}$ have the
following expressions
$$
\begin{array}{lll}
\displaystyle \widetilde{E}_{j}^{k}=\widetilde{\zeta}_{k}
\exp\left(i\beta j h \right),
\;\;
\displaystyle \widetilde{R}_{j}^{k}=\widetilde{\eta}_{k}
\exp\left(i\beta j h \right),
\end{array}
$$
respectively. Substituting the above two expressions into (33) yields
$$
\begin{array}{lll}
\displaystyle \mathcal
{\widetilde{Q}}\widetilde{\zeta}_{k+1}=&\displaystyle \mathcal
{\widetilde{P}}\widetilde{\zeta}_{k} -4\sin^2\left( \frac{\beta
h}{2}\right)\left[1+\frac{1}{3}\sin^2\left( \frac{\beta
h}{2}\right)+\frac{8}{45}\sin^4\left( \frac{\beta
h}{2}\right)\right]\vspace{0.2 cm}\\& \displaystyle\times
\sum\limits_{\ell=2}^{k+1}
g_{\ell}^{(\alpha,\beta)}
\widetilde{\zeta}_{k+1-\ell}+\tau\left[1-\frac{4}{35}\sin^6\left(
\frac{\beta h}{2}\right)\right]\eta_{k+1}.
\end{array}\eqno(34)
$$

{\it{\bf Lemma 12.} Let $\widetilde{\zeta}_{k+1}$
$(k=0,1,\cdots,N-1)$ be the solution of equation (34), under
condition (27), then there exists a positive constant
$\widetilde{C}_2$ such that
$$
\begin{array}{lll}\displaystyle
\left|\widetilde{\zeta}_{k+1}\right|\leq
\widetilde{C}_{2}(k+1)\tau\left|\widetilde{\eta}_1\right|,
\;\;\;k=0,1,\cdots,N-1.
\end{array}
$$}

\begin{proof}
The proof is almost the same as that of Lemma 11, so is omitted here.
\end{proof}

 {\it{\bf Theorem 6.} Under condition (27), the difference scheme (14) is convergent with
order $O(\tau^2+h^8)$.}

\begin{proof} The proof is the same as that of Theorem 5, so is left out here.
\end{proof}

{\it{\bf Remark 2:} In view of conditions (24) and (27), we find
that if $\alpha\in(0,2-\sqrt{2}]$ and $\beta\in(0,2-\sqrt{2}]$, then difference schemes (13)
and (14) are both unconditionally stable. If $\alpha\in(2-\sqrt{2},1)$ or $\beta\in(2-\sqrt{2},1)$, the difference
schemes (13) and (14) are both conditionally stable provided that the stability
conditions are (24) and (27) are still satisfied.}

 \section{Numerical example}
\qquad In this section we list the numerical results of the finite difference schemes in the paper and in \cite{MAD} on one test problem.
 We show the convergence orders and
stability of the methods developed in this paper by performing the mentioned schemes for different values
 of $\alpha,\beta$, $\tau$ and $h$.  All our tests were done in MATLAB.
 The maximum norm error between the exact solution and the numerical
solution is defined as follows:
$$
\begin{array}{lll}
\displaystyle e_{\infty}(\tau,h)=\max_{1\leq j\leq M-1,\;0\leq k\leq N}\left|u_j^k-u(x_j,t_k)\right|.
\end{array}
$$

Define the convergence orders in the temporal direction by
$$
\begin{array}{lll}
\displaystyle \textmd{T-order}=\log_2\left(\frac{e_{\infty}(2\tau,h)}{e_{\infty}(\tau,h)}\right),
\end{array}
$$
and in the spatial direction by
$$
\begin{array}{lll}
\displaystyle \textmd{S-order}=\log_\frac{1}{1-2h}\left(\frac{e_{\infty}(\tau,\frac{1}{1-2h}h)}{e_{\infty}(\tau,h)}\right),
\end{array}
$$
 respectively.

{\bf Example:} Consider the following modified anomalous subdiffusion equation
$$
\begin{array}{lll} \displaystyle
 \frac{\partial{{}u(x,t)}}{\partial{t}}
 =\left(\,_{RL}D_{0,t}^{1-\alpha}+
\,_{RL}D_{0,t}^{1-\beta}
 \right)\left[
 \frac{\partial^2 u(x,t)}{\partial{x^2}}\right]
  +f(x,t),\;\;\;0<x<1,\;\;\;0< t\leq 1$$,
\end{array}
$$
where
$$
\begin{array}{lll} \displaystyle
 f(x,t)&=&\displaystyle(\alpha+\beta+2)t^{\alpha+\beta+1}x^{12}(1-x)^{12}\sin(\pi
 x)+x^{10}(1-x)^{10}\vspace{0.2cm}\\&&\displaystyle
\times\left[\sin(\pi
x)(\pi^2x^2(1-x)^2-552x^2+552x-132)\right.\vspace{0.2cm}\\&&\displaystyle\left.-24\pi
x\cos(\pi
x)(2x^2-3x+1)\right]\vspace{0.2cm}\\
&& \displaystyle\times\left[
\frac{\Gamma(\alpha+\beta+3)}{\Gamma(2\alpha+\beta+2)}t^{2\alpha+\beta+1}
+\frac{\Gamma(\alpha+\beta+3)}{\Gamma(2\beta+\alpha+2)}t^{2\beta+\alpha+1}
 \right].
\end{array}
$$
 Its exact solution is $u(x, t)
=t^{\alpha+\beta+2}x^{12}(1-x)^{12}\sin(\pi x) $, which satisfies the
initial and boundary value conditions. This equation for
describing processes that become less anomalous as time progresses by the inclusion of
a second fractional time derivative acting on the diffusion term. The subdiffusive motion is characterized by an asymptotic longtime
behavior of the mean square displacement of the form \cite{L}
$$
\begin{array}{lll}
\displaystyle \langle x^2(t)\rangle\;\sim\;\frac{2}{\Gamma(1+\alpha)}t^{\alpha}+\frac{2}{\Gamma(1+\beta)}t^{\beta},\;t\longrightarrow\infty.
\end{array}
$$
A possible application of this equation is in econophysics.
 In particular the crossover between more and less anomalous behavior has been observed in
the volatility of some share prices \cite{MMPW}.

Here, we compare the numerical results of the finite difference schemes (13), (14)
with those of the numerical scheme in \cite{MAD}. The maximum-norm error, temporal and spatial
convergence orders, and CPU time for these finite difference schemes are listed in
Tables 1--4 for different $\alpha$, $\beta$. From these tables,
it is clear to see that the finite difference schemes (13) and (14) provide much more accuracy
and do not lead to additional computational requirements than those in \cite{MAD}
for the same grid sizes. Furthermore, one can seen that the
computational convergence orders are close to theoretical convergence orders, i.e., the convergence orders of the finite difference schemes
 (13) and (14) in temporal direction are both second-orders, in spatial direction are sixth-order and eight-order, respectively.

\begin{figure}[!htbp]
\centering
 \includegraphics[width=16 cm]{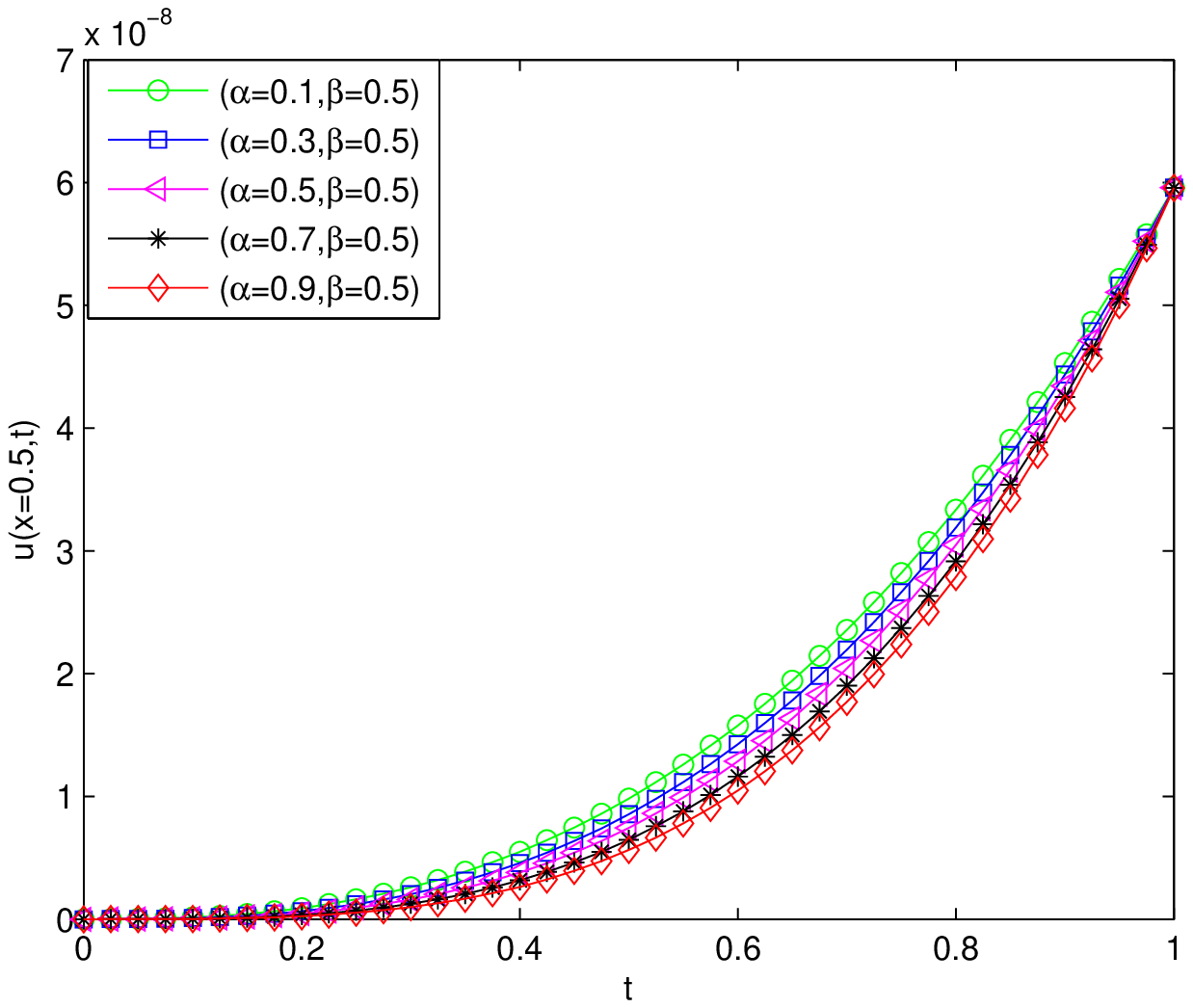}\\
  \caption{The numerical solution behaviours at $x=0.5$ by difference scheme (13) for different $\alpha$
  and fixed $\beta=0.5$ with $\tau=\frac{1}{40}$ and $h=\frac{1}{100}$.}
  \label{fig.1}
\end{figure}

\begin{figure}[!htbp]
\centering
 \includegraphics[width=16 cm]{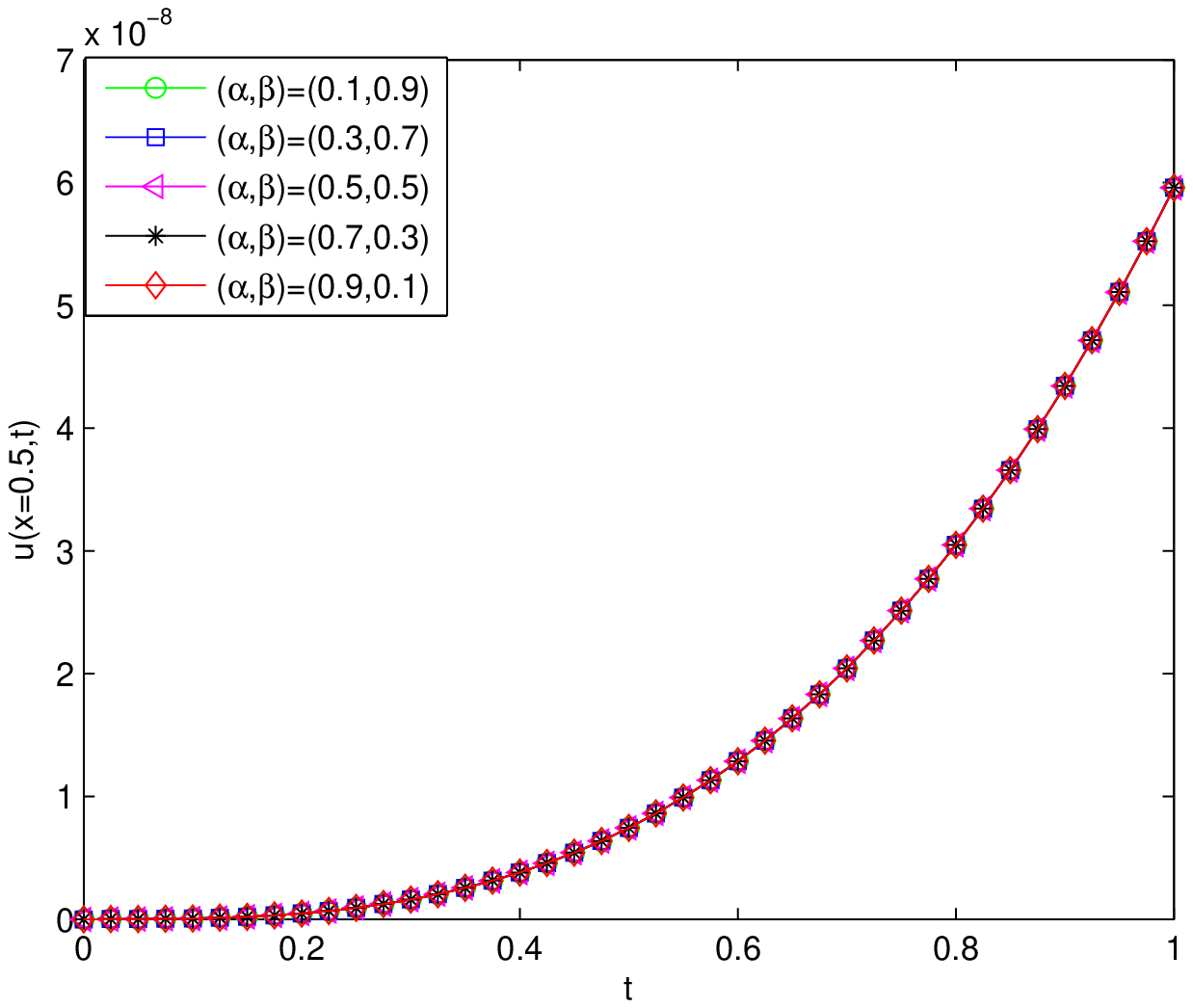}\\
  \caption{The numerical solution behaviours at $x=0.5$ by difference scheme (13) for different pairs $(\alpha,\beta)$
  (which satisfy $\alpha+\beta=1$) with $\tau=\frac{1}{40}$ and $h=\frac{1}{100}$.}
  \label{fig.2}
\end{figure}

\begin{figure}[!htbp]
\centering
 \includegraphics[width=16 cm]{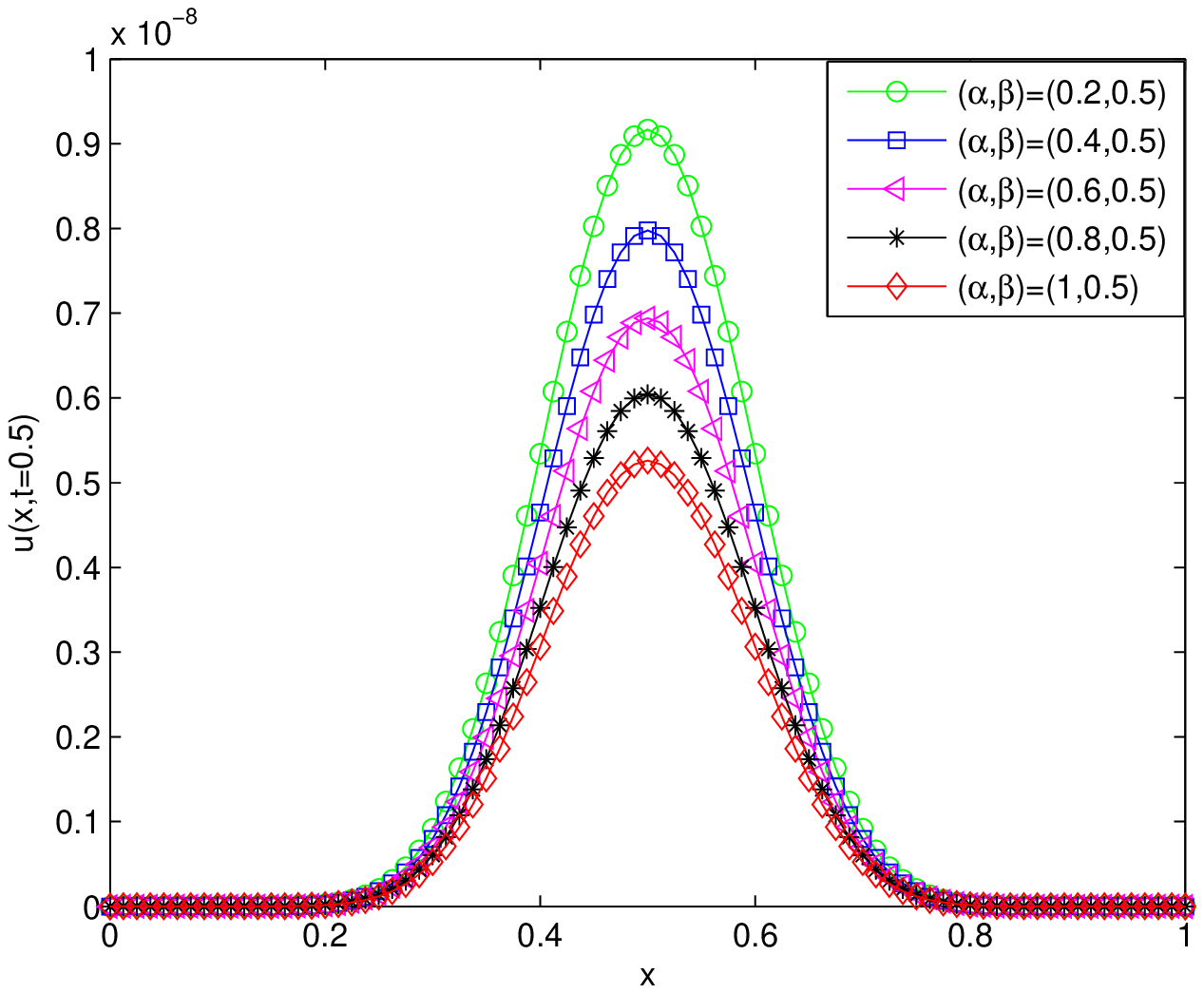}\\
  \caption{The numerical solution behaviours at $t=0.5$ by difference scheme (13) for different $\alpha$
  and fixed $\beta=0.5$ with $\tau=\frac{1}{50}$ and $h=\frac{1}{80}$.}
  \label{fig.3}
\end{figure}

\begin{figure}[!htbp]
\centering
 \includegraphics[width=16 cm]{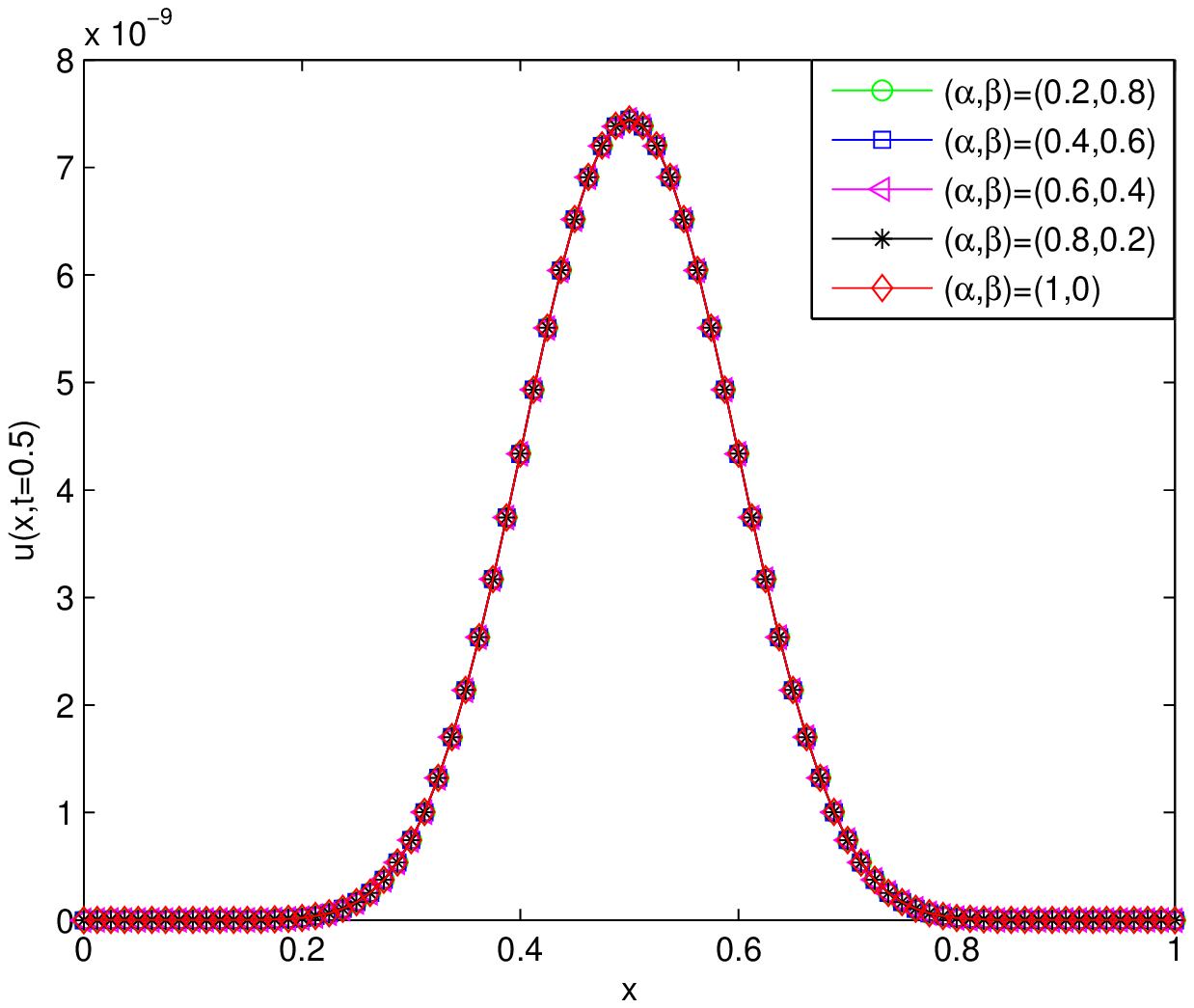}\\
  \caption{The numerical solution behaviours at $t=0.5$ by difference scheme (13) for different pairs $(\alpha,\beta)$
  (which satisfy $\alpha+\beta=1$) with $\tau=\frac{1}{40}$ and $h=\frac{1}{100}$.}
  \label{fig.4}
\end{figure}

Next, we display the numerical solutions profiles for difference cases by Figs \ref{fig.1}--\ref{fig.4}.
From these Figs, it is clear that the equation in the paper exhibits anomalous diffusion behaviours and the
fractional differential equations are characterised by a heavy tail (see Figs \ref{fig.3} and \ref{fig.4}).
For the probability density
function associated with such diffusion process is no longer Gaussian
but is replaced by a more general L\'{e}vy distribution. A distribution
which can exhibit heavy tails with a power law decay as
opposed to the thin exponentially decaying tails of a Gaussian distribution \cite{MK} and
 resulting in long-range dependence.
In addition, We found an interesting phenomenon is that these numerical solutions for different pairs $(\alpha,\beta)$
 show almost the same behaviors as long as meet the condition $\alpha+\beta=1$ (see Figs \ref{fig.2} and \ref{fig.4}).

\red{
\begin{table}[!htbp]\renewcommand\arraystretch{1.4}
 \begin{center}
 \caption{ The comparison of the difference scheme (13) with difference scheme in \cite{MAD}
for $h=1/1000$.}\vspace{1
cm}
 \begin{footnotesize}
\begin{tabular}{c c c c c @{}p{2mm}@{} c c c}\hline\\
\multicolumn{1}{c}{\multirow{2}{*}{$(\alpha,\beta)$}}   & \multicolumn{1}{c}{\multirow{2}{*}{$\tau$}}&
 \multicolumn{3}{c} {Finite difference scheme (13)}&& \multicolumn{3}{c} {Finite difference scheme in \cite{MAD}}\\
 \cline{3-5}\cline{7-9}
 \multicolumn{1}{c}{}&\multicolumn{1}{c}{}&\multicolumn{1}{c}{\textrm{$e_{\infty}(\tau,h)$}}&$\textmd{T-order}$ &CPU time (s)&&\multicolumn{1}{c}{\textrm{$e_{\infty}(\tau,h)$}}&$\textmd{T-order}$ &CPU time (s)\\
 \hline
 \multicolumn{1}{c}{(0.25,0.15)}&\multicolumn{1}{c}{$\frac{1}{4}$}&\multicolumn{1}{c}{9.1447e-010}&---&2.270
 &&\multicolumn{1}{c}{ 1.9599e-008}&--- &3.885\\
  \cline{2-9}
\multicolumn{1}{c}{}&\multicolumn{1}{c}{$\frac{1}{8}$}&\multicolumn{1}{c}{ 2.3336e-010}& 1.9704&2.849
 &&\multicolumn{1}{c}{1.0226e-008}&0.9385 &8.054\\
  \cline{2-9}
  \multicolumn{1}{c}{}&\multicolumn{1}{c}{$\frac{1}{16}$}&\multicolumn{1}{c}{ 5.9079e-011}&   1.9818&10.851
 &&\multicolumn{1}{c}{ 5.2218e-009}&0.9696&18.776\\
  \cline{2-9}
  \multicolumn{1}{c}{}&\multicolumn{1}{c}{$\frac{1}{32}$}&\multicolumn{1}{c}{1.4886e-011}& 1.9887& 10.278
 &&\multicolumn{1}{c}{ 2.6386e-009}&0.9848 &47.736\\
  \hline
  \multicolumn{1}{c}{(0.25,0.35)}&\multicolumn{1}{c}{$\frac{1}{4}$}&\multicolumn{1}{c}{ 1.3779e-009}&---&2.092
 &&\multicolumn{1}{c}{2.2275e-008}&--- &3.762\\
  \cline{2-9}
\multicolumn{1}{c}{}&\multicolumn{1}{c}{$\frac{1}{8}$}&\multicolumn{1}{c}{ 3.4896e-010}&1.9813&2.814
 &&\multicolumn{1}{c}{1.1805e-008}&0.9160 &8.083\\
  \cline{2-9}
  \multicolumn{1}{c}{}&\multicolumn{1}{c}{$\frac{1}{16}$}&\multicolumn{1}{c}{ 8.7830e-011}&  1.9903&4.673
 &&\multicolumn{1}{c}{6.0722e-009}& 0.9591&18.489\\
  \cline{2-9}
  \multicolumn{1}{c}{}&\multicolumn{1}{c}{$\frac{1}{32}$}&\multicolumn{1}{c}{2.2034e-011}& 1.9950&10.323
 &&\multicolumn{1}{c}{3.0790e-009}&  0.9798 &46.969\\
  \hline
  \multicolumn{1}{c}{(0.25,0.55)}&\multicolumn{1}{c}{$\frac{1}{4}$}&\multicolumn{1}{c}{1.8262e-009}&---&2.102
 &&\multicolumn{1}{c}{2.4748e-008}&--- &3.742\\
  \cline{2-9}
\multicolumn{1}{c}{}&\multicolumn{1}{c}{$\frac{1}{8}$}&\multicolumn{1}{c}{4.6262e-010}& 1.9809&2.804
 &&\multicolumn{1}{c}{ 1.3307e-008}& 0.8951&8.130\\
  \cline{2-9}
  \multicolumn{1}{c}{}&\multicolumn{1}{c}{$\frac{1}{16}$}&\multicolumn{1}{c}{ 1.1644e-010}& 1.9902&4.667
 &&\multicolumn{1}{c}{6.8935e-009}& 0.9489 &18.731\\
  \cline{2-9}
  \multicolumn{1}{c}{}&\multicolumn{1}{c}{$\frac{1}{32}$}&\multicolumn{1}{c}{ 2.9210e-011}& 1.9951&10.214
 &&\multicolumn{1}{c}{3.5075e-009}& 0.9748 &48.639\\
  \hline
\end{tabular}
 \end{footnotesize}
 \end{center}
 \end{table}}

\red{
\begin{table}[!htbp]\renewcommand\arraystretch{1.4}
 \begin{center}
 \caption{ The comparison of the difference scheme (13) with difference scheme in \cite{MAD}
for $\tau=1/200$.}\vspace{1
cm}
 \begin{footnotesize}
\begin{tabular}{c c c c c @{}p{2mm}@{} c c c}\hline\\
\multicolumn{1}{c}{\multirow{2}{*}{$(\alpha,\beta)$}}   & \multicolumn{1}{c}{\multirow{2}{*}{$h$}}&
 \multicolumn{3}{c} {Finite difference scheme (13)}&& \multicolumn{3}{c} {Finite difference scheme in \cite{MAD}}\\
 \cline{3-5}\cline{7-9}
 \multicolumn{1}{c}{}&\multicolumn{1}{c}{}&\multicolumn{1}{c}{\textrm{$e_{\infty}(\tau,h)$}}&$\textmd{S-order}$ &CPU time (s)&&\multicolumn{1}{c}{\textrm{$e_{\infty}(\tau,h)$}}&$\textmd{S-order}$ &CPU time (s)\\
 \hline
 \multicolumn{1}{c}{(0.4,0.1)}&\multicolumn{1}{c}{$\frac{1}{12}$}&\multicolumn{1}{c}{ 1.8047e-010}&---&3.398
 &&\multicolumn{1}{c}{ 3.7736e-010}&--- &3.802\\
  \cline{2-9}
\multicolumn{1}{c}{}&\multicolumn{1}{c}{$\frac{1}{14}$}&\multicolumn{1}{c}{ 7.5031e-011}&  5.6935&3.352
 &&\multicolumn{1}{c}{ 3.2936e-010}&  0.8826 &4.059\\
  \cline{2-9}
  \multicolumn{1}{c}{}&\multicolumn{1}{c}{$\frac{1}{16}$}&\multicolumn{1}{c}{3.4678e-011}&   5.7799&3.529
 &&\multicolumn{1}{c}{3.3801e-010}&not convergent&4.364\\
  \cline{2-9}
  \multicolumn{1}{c}{}&\multicolumn{1}{c}{$\frac{1}{18}$}&\multicolumn{1}{c}{1.7276e-011}& 5.9159&  3.584
 &&\multicolumn{1}{c}{  3.7944e-010}&not convergent&4.643\\
  \hline
  \multicolumn{1}{c}{(0.4,0.3)}&\multicolumn{1}{c}{$\frac{1}{12}$}&\multicolumn{1}{c}{1.8011e-010}&---&3.293
 &&\multicolumn{1}{c}{4.3226e-010}&--- &3.825\\
  \cline{2-9}
\multicolumn{1}{c}{}&\multicolumn{1}{c}{$\frac{1}{14}$}&\multicolumn{1}{c}{7.4697e-011}&5.7095&3.381
 &&\multicolumn{1}{c}{ 3.9031e-010}& 0.6622&4.102\\
  \cline{2-9}
  \multicolumn{1}{c}{}&\multicolumn{1}{c}{$\frac{1}{16}$}&\multicolumn{1}{c}{ 3.4353e-011}& 5.8170& 3.432
 &&\multicolumn{1}{c}{ 4.1165e-010}& not convergent&4.354\\
  \cline{2-9}
  \multicolumn{1}{c}{}&\multicolumn{1}{c}{$\frac{1}{18}$}&\multicolumn{1}{c}{1.6955e-011}& 5.9951&3.555
 &&\multicolumn{1}{c}{ 4.6002e-010}&not convergent &4.738\\
  \hline
  \multicolumn{1}{c}{(0.4,0.5)}&\multicolumn{1}{c}{$\frac{1}{12}$}&\multicolumn{1}{c}{1.7987e-010}&---&3.392
 &&\multicolumn{1}{c}{4.8705e-010}&--- &3.872\\
  \cline{2-9}
\multicolumn{1}{c}{}&\multicolumn{1}{c}{$\frac{1}{14}$}&\multicolumn{1}{c}{7.4486e-011}& 5.7192&3.368
 &&\multicolumn{1}{c}{4.5116e-010}&  0.4966&4.090\\
  \cline{2-9}
  \multicolumn{1}{c}{}&\multicolumn{1}{c}{$\frac{1}{16}$}&\multicolumn{1}{c}{3.4153e-011}& 5.8395&3.458
 &&\multicolumn{1}{c}{ 4.9230e-010}& not convergent &4.327\\
  \cline{2-9}
  \multicolumn{1}{c}{}&\multicolumn{1}{c}{$\frac{1}{18}$}&\multicolumn{1}{c}{ 1.6760e-011}& 6.0438&3.543
 &&\multicolumn{1}{c}{ 5.4059e-010}& not convergent& 4.660\\
  \hline
\end{tabular}
 \end{footnotesize}
 \end{center}
 \end{table}}

\red{
\begin{table}[!htbp]\renewcommand\arraystretch{1.4}
 \begin{center}
 \caption{ The comparison of the difference scheme (14) with difference scheme in \cite{MAD}
for $h=1/500$.}\vspace{1
cm}
 \begin{footnotesize}
\begin{tabular}{c c c c c @{}p{2mm}@{} c c c}\hline\\
\multicolumn{1}{c}{\multirow{2}{*}{$(\alpha,\beta)$}}   & \multicolumn{1}{c}{\multirow{2}{*}{$\tau$}}&
 \multicolumn{3}{c} {Finite difference scheme (14)}&& \multicolumn{3}{c} {Finite difference scheme in \cite{MAD}}\\
 \cline{3-5}\cline{7-9}
 \multicolumn{1}{c}{}&\multicolumn{1}{c}{}&\multicolumn{1}{c}{\textrm{$e_{\infty}(\tau,h)$}}&$\textmd{T-order}$ &CPU time (s)&&\multicolumn{1}{c}{\textrm{$e_{\infty}(\tau,h)$}}&$\textmd{T-order}$ &CPU time (s)\\
 \hline
 \multicolumn{1}{c}{(0.45,0.15)}&\multicolumn{1}{c}{$\frac{1}{4}$}&\multicolumn{1}{c}{1.3338e-009}&---&0.508
 &&\multicolumn{1}{c}{2.2051e-008}&--- &0.654\\
  \cline{2-9}
\multicolumn{1}{c}{}&\multicolumn{1}{c}{$\frac{1}{8}$}&\multicolumn{1}{c}{ 3.3783e-010}& 1.9812&0.747
 &&\multicolumn{1}{c}{ 1.1676e-008}&0.9173 & 1.486\\
  \cline{2-9}
  \multicolumn{1}{c}{}&\multicolumn{1}{c}{$\frac{1}{16}$}&\multicolumn{1}{c}{8.5051e-011}& 1.9899& 1.397
 &&\multicolumn{1}{c}{ 6.0042e-009}& 0.9595&3.484\\
  \cline{2-9}
  \multicolumn{1}{c}{}&\multicolumn{1}{c}{$\frac{1}{32}$}&\multicolumn{1}{c}{2.1342e-011}& 1.9946&3.093
 &&\multicolumn{1}{c}{  3.0440e-009}& 0.9800 &9.355\\
  \hline
  \multicolumn{1}{c}{(0.45,0.35)}&\multicolumn{1}{c}{$\frac{1}{4}$}&\multicolumn{1}{c}{1.8801e-009}&---&0.468
 &&\multicolumn{1}{c}{ 2.5000e-008}&--- &0.661\\
  \cline{2-9}
\multicolumn{1}{c}{}&\multicolumn{1}{c}{$\frac{1}{8}$}&\multicolumn{1}{c}{ 4.7658e-010}&1.9800&0.729
 &&\multicolumn{1}{c}{ 1.3456e-008}& 0.8937&1.430\\
  \cline{2-9}
  \multicolumn{1}{c}{}&\multicolumn{1}{c}{$\frac{1}{16}$}&\multicolumn{1}{c}{1.2000e-010}& 1.9897&1.389
 &&\multicolumn{1}{c}{ 6.9737e-009}&  0.9483&3.427\\
  \cline{2-9}
  \multicolumn{1}{c}{}&\multicolumn{1}{c}{$\frac{1}{32}$}&\multicolumn{1}{c}{3.0108e-011}& 1.9948&3.132
 &&\multicolumn{1}{c}{3.5491e-009}&  0.9745&9.412\\
  \hline
  \multicolumn{1}{c}{(0.45,0.55)}&\multicolumn{1}{c}{$\frac{1}{4}$}&\multicolumn{1}{c}{2.4266e-009}&---&0.464
 &&\multicolumn{1}{c}{ 2.7710e-008}&--- &0.649\\
  \cline{2-9}
\multicolumn{1}{c}{}&\multicolumn{1}{c}{$\frac{1}{8}$}&\multicolumn{1}{c}{6.1660e-010}& 1.9765&0.743
 &&\multicolumn{1}{c}{ 1.5161e-008}&    0.8700&1.419\\
  \cline{2-9}
  \multicolumn{1}{c}{}&\multicolumn{1}{c}{$\frac{1}{16}$}&\multicolumn{1}{c}{1.5543e-010}&1.9881&1.399
 &&\multicolumn{1}{c}{7.9204e-009}&0.9367 &3.511\\
  \cline{2-9}
  \multicolumn{1}{c}{}&\multicolumn{1}{c}{$\frac{1}{32}$}&\multicolumn{1}{c}{3.9020e-011}&1.9940&3.121
 &&\multicolumn{1}{c}{4.0469e-009}&0.9688 &9.789\\
  \hline
\end{tabular}
 \end{footnotesize}
 \end{center}
 \end{table}}

\red{
\begin{table}[!htbp]\renewcommand\arraystretch{1.4}
 \begin{center}
 \caption{ The comparison of the difference scheme (14) with difference scheme in \cite{MAD}
for $\tau=1/160$.}\vspace{1
cm}
 \begin{footnotesize}
\begin{tabular}{c c c c c @{}p{2mm}@{} c c c}\hline\\
\multicolumn{1}{c}{\multirow{2}{*}{$(\alpha,\beta)$}}   & \multicolumn{1}{c}{\multirow{2}{*}{$h$}}&
 \multicolumn{3}{c} {Finite difference scheme (14)}&& \multicolumn{3}{c} {Finite difference scheme in \cite{MAD}}\\
 \cline{3-5}\cline{7-9}
 \multicolumn{1}{c}{}&\multicolumn{1}{c}{}&\multicolumn{1}{c}{\textrm{$e_{\infty}(\tau,h)$}}&$\textmd{S-order}$ &CPU time (s)&&\multicolumn{1}{c}{\textrm{$e_{\infty}(\tau,h)$}}&$\textmd{S-order}$ &CPU time (s)\\
 \hline
 \multicolumn{1}{c}{(0.2,0.1)}&\multicolumn{1}{c}{$\frac{1}{14}$}&\multicolumn{1}{c}{3.1090e-011}&---&5.160
 &&\multicolumn{1}{c}{3.5451e-010}&--- &2.700\\
  \cline{2-9}
\multicolumn{1}{c}{}&\multicolumn{1}{c}{$\frac{1}{16}$}&\multicolumn{1}{c}{1.1703e-011}& 7.3169&2.295
 &&\multicolumn{1}{c}{3.6480e-010}& not convergent &2.850\\
  \cline{2-9}
  \multicolumn{1}{c}{}&\multicolumn{1}{c}{$\frac{1}{18}$}&\multicolumn{1}{c}{4.6941e-012}&7.7561&2.368
 &&\multicolumn{1}{c}{ 4.1230e-010}&not convergent&3.009\\
  \cline{2-9}
  \multicolumn{1}{c}{}&\multicolumn{1}{c}{$\frac{1}{20}$}&\multicolumn{1}{c}{2.0284e-012}& 7.9637& 2.433
 &&\multicolumn{1}{c}{  4.3941e-010}&not convergent&3.293\\
  \hline
  \multicolumn{1}{c}{(0.2,0.3)}&\multicolumn{1}{c}{$\frac{1}{14}$}&\multicolumn{1}{c}{3.0823e-011}&---&2.187
 &&\multicolumn{1}{c}{ 4.2057e-010}&--- &2.639\\
  \cline{2-9}
\multicolumn{1}{c}{}&\multicolumn{1}{c}{$\frac{1}{16}$}&\multicolumn{1}{c}{1.1437e-011}&7.4245&2.245
 &&\multicolumn{1}{c}{ 4.5120e-010}&not convergent&2.781\\
  \cline{2-9}
  \multicolumn{1}{c}{}&\multicolumn{1}{c}{$\frac{1}{18}$}&\multicolumn{1}{c}{4.6804e-012}&7.5857&2.311
 &&\multicolumn{1}{c}{ 4.9954e-010}& not convergent& 3.049\\
  \cline{2-9}
  \multicolumn{1}{c}{}&\multicolumn{1}{c}{$\frac{1}{20}$}&\multicolumn{1}{c}{2.1812e-012}&7.2466& 2.360
 &&\multicolumn{1}{c}{ 5.2660e-010}&not convergent & 3.220\\
  \hline
  \multicolumn{1}{c}{(0.2,0.5)}&\multicolumn{1}{c}{$\frac{1}{14}$}&\multicolumn{1}{c}{3.0548e-011}&---&2.171
 &&\multicolumn{1}{c}{ 4.8574e-010}&--- &2.646\\
  \cline{2-9}
\multicolumn{1}{c}{}&\multicolumn{1}{c}{$\frac{1}{16}$}&\multicolumn{1}{c}{1.1164e-011}& 7.5383&2.241
 &&\multicolumn{1}{c}{ 5.3745e-010}& not convergent&2.805\\
  \cline{2-9}
  \multicolumn{1}{c}{}&\multicolumn{1}{c}{$\frac{1}{18}$}&\multicolumn{1}{c}{4.8168e-012}&7.1367&2.305
 &&\multicolumn{1}{c}{  5.8571e-010}& not convergent &3.014\\
  \cline{2-9}
  \multicolumn{1}{c}{}&\multicolumn{1}{c}{$\frac{1}{20}$}&\multicolumn{1}{c}{2.3377e-012}&6.8616&2.366
 &&\multicolumn{1}{c}{  6.1272e-010}& not convergent& 3.271\\
  \hline
\end{tabular}
 \end{footnotesize}
 \end{center}
 \end{table}}

\section{Conclusion}
\quad\ In this paper, we establish two high-order compact finite
difference schemes for the modified anomalous subdiffusion equation.
The stability and convergence conditions of the difference schemes
 are given by using the Fourier method. Finally, numerical experiments have
been carried out to support the theoretical claims. These methods
and techniques can be extended in a straightforward way to two
or three spatial dimensional cases.\\

{\it\bf Appendix I:}

The forms of the matrices $A,B,$ $\widetilde{A}, \widetilde{B}$, $C_k$ and $\widetilde{C_k}$ are list as follows:

$$ \displaystyle
 A= \left(
  \begin{array}{ccccccccc}
    \frac{14}{15} & \frac{2}{45} &  -\frac{1}{90}& 0 & \cdots& 0 &
    0& 0 & 0\vspace{0.2 cm}\\
  \frac{2}{45}& \frac{14}{15}& \frac{2}{45}
     &  -\frac{1}{90} & 0& \cdots &0 &
      0 & 0\vspace{0.2 cm}\\
  -\frac{1}{90} & \frac{2}{45}  & \frac{14}{15}
     &\frac{2}{45} & -\frac{1}{90} & 0 & \cdots&
     0 & 0\vspace{0.2 cm}\\
   0& -\frac{1}{90} & \frac{2}{45}
     &\frac{14}{15} & \frac{2}{45} & -\frac{1}{90} & 0&   \cdots & 0\vspace{0.4 cm}\\
     \vdots & \ddots & \ddots & \ddots &\ddots &\ddots & \ddots & \ddots &\vdots\vspace{0.4 cm}\\
  0&\cdots&   0& -\frac{1}{90} & \frac{2}{45} &
    \frac{14}{15}  & \frac{2}{45} & -\frac{1}{90} & 0\vspace{0.2 cm}\\
    0& 0 &\cdots& 0 &  -\frac{1}{90}
    & \frac{2}{45} & \frac{14}{15} & \frac{2}{45} &-\frac{1}{90}\vspace{0.2 cm}\\
   0&0& 0   & \cdots& 0 & -\frac{1}{90} & \frac{2}{45} & \frac{14}{15}  &  \frac{2}{45}
     \vspace{0.2 cm}\\
   0&  0 & 0 & 0
    & \cdots& 0&  -\frac{1}{90}& \frac{2}{45}
    &  \frac{14}{15}  \vspace{0.2 cm}\\
  \end{array}
\right),
$$
$$ \displaystyle
B= \left(
  \begin{array}{ccccccccc}
     -\frac{5}{2} & \frac{4}{3} &  -\frac{1}{12}& 0 & \cdots& 0 &
    0& 0& 0 \vspace{0.2 cm}\\
  \frac{4}{3}&  -\frac{5}{2}& \frac{4}{3}
     &  -\frac{1}{12} & 0& \cdots &0 &
      0 & 0\vspace{0.2 cm}\\
  -\frac{1}{12} & \frac{4}{3}  &  -\frac{5}{2}
     &\frac{4}{3} & -\frac{1}{12} & 0 & \cdots&
     0 & 0\vspace{0.2 cm}\\
   0& -\frac{1}{12} & \frac{4}{3}
     & -\frac{5}{2} & \frac{4}{3} & -\frac{1}{12} & 0&  \cdots &0\vspace{0.4 cm}\\
    \vdots & \ddots & \ddots & \ddots & \ddots &\ddots &\ddots & \ddots & \vdots\vspace{0.4 cm}\\
    0 &\cdots& 0 & -\frac{1}{12} & \frac{4}{3} &
     -\frac{5}{2}  & \frac{4}{3} & -\frac{1}{12}&
     0\vspace{0.2 cm}\\
   0& 0 &  \cdots
    & 0 & -\frac{1}{12} & \frac{4}{3} &
     -\frac{5}{2}  & \frac{4}{3} & -\frac{1}{12}\vspace{0.2 cm}\\
    0& 0& 0 &\cdots& 0 &-\frac{1}{12} & \frac{4}{3} &
     -\frac{5}{2}  & \frac{4}{3} \vspace{0.2 cm}\\
    0& 0& 0 & 0
    & \cdots& 0&  -\frac{1}{12}& \frac{4}{3}
    &   -\frac{5}{2}   \vspace{0.2 cm}\\
  \end{array}
\right),
$$

$$ \displaystyle
\widetilde{ A}= \left(
  \begin{array}{ccccccccccc}
    \frac{27}{28} & \frac{3}{112}&  -\frac{3}{280}&
     \frac{1}{560}& 0& \cdots &0& 0&  0 & 0&  0\vspace{0.2 cm}\\
   \frac{3}{112}& \frac{27}{28}&
   \frac{3}{112} & -\frac{3}{280}& \frac{1}{560} &0 &\cdots &0& 0& 0&  0\vspace{0.2 cm}\\
     -\frac{3}{280}  & \frac{3}{112}
     &\frac{27}{28} & \frac{3}{112} & -\frac{3}{280} & \frac{1}{560}
     &0 &
    \cdots& 0& 0&  0\vspace{0.2 cm}\\
   \frac{1}{560} &-\frac{3}{280}
     & \frac{3}{112} & \frac{27}{28}& \frac{3}{112} & -\frac{3}{280}&  \frac{1}{560}&
    0& \cdots & 0&  0\vspace{0.2 cm}\\
  0 & \frac{1}{560}
     &-\frac{3}{280} & \frac{3}{112} & \frac{27}{28} & \frac{3}{112}&  -\frac{3}{280}&
    \frac{1}{560}& 0& 0&  0\vspace{0.4 cm}\\
    \vdots & \ddots & \ddots & \ddots & \ddots &\ddots &\ddots & \ddots & \ddots&
    \ddots& \vdots\vspace{0.4 cm}\\
    0&\cdots&0 & \frac{1}{560} &  -\frac{3}{280}& \frac{3}{112}& \frac{27}{28}&
    \frac{3}{112} & -\frac{3}{280}& \frac{1}{560} &
    0\vspace{0.2 cm}\\
    0& 0&\cdots &0& \frac{1}{560} &  -\frac{3}{280}& \frac{3}{112}& \frac{27}{28}&
     \frac{3}{112} & -\frac{3}{280}& \frac{1}{560}\vspace{0.2 cm}\\
   0& 0 & 0 & \cdots& 0 & \frac{1}{560} &  -\frac{3}{280}& \frac{3}{112}& \frac{27}{28}&
     \frac{3}{112} & -\frac{3}{280}\vspace{0.2 cm}\\
    0& 0 & 0& 0& \cdots& 0 & \frac{1}{560} &  -\frac{3}{280}& \frac{3}{112}& \frac{27}{28}&
     \frac{3}{112} \vspace{0.2 cm}\\
    0& 0 & 0& 0 & 0&\cdots
    & 0&\frac{1}{560} &  -\frac{3}{280}& \frac{3}{112}& \frac{27}{28} \vspace{0.2 cm}\\
  \end{array}
\right),
$$

$$ \displaystyle
\widetilde{ B}= \left(
  \begin{array}{ccccccccccc}
    -\frac{49}{18} & \frac{3}{2}&  -\frac{3}{20}&
     \frac{1}{90}& 0& \cdots &0& 0&  0 & 0&  0\vspace{0.2 cm}\\
  \frac{3}{2}& -\frac{49}{18}&
   \frac{3}{2} & -\frac{3}{20}& \frac{1}{90} &0 &\cdots &0& 0& 0&  0\vspace{0.2 cm}\\
    -\frac{3}{20}  & \frac{3}{2}
     &-\frac{49}{18} & \frac{3}{2} & -\frac{3}{20} & \frac{1}{90}
     &0 &
    \cdots& 0& 0&  0\vspace{0.2 cm}\\
 \frac{1}{90} & -\frac{3}{20}
     & \frac{3}{2} & -\frac{49}{18}& \frac{3}{2} & -\frac{3}{20}&  \frac{1}{90}&
    0& \cdots & 0&  0\vspace{0.2 cm}\\
   0 & \frac{1}{90}
     &-\frac{3}{20} & \frac{3}{2} & -\frac{49}{18} & \frac{3}{2}&  -\frac{3}{20}&
    \frac{1}{90}& 0& \cdots&  0\vspace{0.2 cm}\\
    \vdots & \ddots & \ddots & \ddots & \ddots &\ddots &\ddots & \ddots & \ddots&
    \ddots& \vdots\vspace{0.4 cm}\\
   0& \cdots& 0 & \frac{1}{90} &  -\frac{3}{20}& \frac{3}{2}& -\frac{49}{18}&
    \frac{3}{2} & -\frac{3}{20}& \frac{1}{90} &
    0\vspace{0.2 cm}\\
   0& 0 & \cdots &0& \frac{1}{90} &  -\frac{3}{20}& \frac{3}{2}& -\frac{49}{18}&
     \frac{3}{2} & -\frac{3}{20}& \frac{1}{90}\vspace{0.2 cm}\\
  0& 0 & 0 & \cdots& 0 & \frac{1}{90} &  -\frac{3}{20}& \frac{3}{2}& -\frac{49}{18}&
     \frac{3}{2} & -\frac{3}{20}\vspace{0.2 cm}\\
   0& 0 & 0& 0& \cdots& 0 & \frac{1}{90} &  -\frac{3}{20}& \frac{3}{2}& -\frac{49}{18}&
     \frac{3}{2} \vspace{0.2 cm}\\
    0& 0 &0& 0 & 0&\cdots
    & 0&\frac{1}{90} &  -\frac{3}{20}& \frac{3}{2}& -\frac{49}{18} \vspace{0.2 cm}\\
  \end{array}
\right),
$$

$$
\begin{array}{l}\displaystyle
C_k=\left(
      \begin{array}{c}
       \left(\frac{1}{90}-\frac{1}{12}g_{0}^{(\alpha,\beta)}\right)u_{-1}^{k+1}
   -\left(\frac{1}{90}+\frac{1}{12}g_{1}^{(\alpha,\beta)}\right)
   u_{-1}^{k}
   -\frac{1}{12}\sum\limits_{\ell=2}^{k+1}
g_{\ell}^{(\alpha,\beta)}u_{-1}^{k+1-\ell}\vspace{0.2 cm}\\-
\left(\frac{2}{45}-\frac{4}{3}g_{0}^{(\alpha,\beta)}\right)u_{0}^{k+1}
+
\left(\frac{2}{45}+\frac{4}{3}g_{1}^{(\alpha,\beta)}\right)u_{0}^{k}
   +\frac{4}{3}\sum\limits_{\ell=2}^{k+1}
g_{\ell}^{(\alpha,\beta)}u_{0}^{k+1-\ell}
\vspace{0.2 cm}\\
-\frac{1}{90}\tau f_{-1}^{k+\frac{1}{2}}+\frac{2}{45}\tau
f_{0}^{k+\frac{1}{2}},
        \vspace{0.4 cm}\\
      \left(\frac{1}{90}-\frac{1}{12}g_{0}^{(\alpha,\beta)}\right)u_{0}^{k+1}
   -\left(\frac{1}{90}+\frac{1}{12}g_{1}^{(\alpha,\beta)}\right)
   u_{0}^{k}
   -\frac{1}{12}\sum\limits_{\ell=2}^{k+1}
g_{\ell}^{(\alpha,\beta)}u_{0}^{k+1-\ell}
-\frac{1}{90}\tau f_{0}^{k+\frac{1}{2}},
\vspace{0.4 cm}\\
0\vspace{0.4 cm}\\
      \vdots\vspace{0.4 cm}\\
       0 \vspace{0.4 cm} \\
            \left(\frac{1}{90}-\frac{1}{12}g_{0}^{(\alpha,\beta)}
            \right)u_{M}^{k+1}
   -\left(\frac{1}{90}+\frac{1}{12}g_{1}^{(\alpha,\beta)}\right)
   u_{M}^{k}
   -\frac{1}{12}\sum\limits_{\ell=2}^{k+1}
g_{\ell}^{(\alpha,\beta)}u_{M}^{k+1-\ell}
-\frac{1}{90}\tau f_{M}^{k+\frac{1}{2}},
\vspace{0.4 cm}\\
  \left(\frac{1}{90}-\frac{1}{12}g_{0}^{(\alpha,\beta)}\right)u_{M+1}^{k+1}
   -\left(\frac{1}{90}+\frac{1}{12}g_{1}^{(\alpha,\beta)}\right)
   u_{M+1}^{k}
   -\frac{1}{12}\sum\limits_{\ell=2}^{k+1}
g_{\ell}^{(\alpha,\beta)}u_{M+1}^{k+1-\ell}\vspace{0.2 cm}\\-
\left(\frac{2}{45}-\frac{4}{3}g_{0}^{(\alpha,\beta)}\right)u_{M}^{k+1}
+
\left(\frac{2}{45}+\frac{4}{3}g_{1}^{(\alpha,\beta)}\right)u_{M}^{k}
   +\frac{4}{3}\sum\limits_{\ell=2}^{k+1}
g_{\ell}^{(\alpha,\beta)}u_{M}^{k+1-\ell}
\vspace{0.2 cm}\\
-\frac{1}{90}\tau f_{M+1}^{k+\frac{1}{2}}+\frac{2}{45}\tau
f_{M}^{k+\frac{1}{2}}
      \end{array}
    \right),
\end{array}
$$

$$
\begin{array}{l}\displaystyle
\widetilde{C}_k=\left(
      \begin{array}{c}
-\left(\frac{1}{560}-\frac{1}{90}g_{0}^{(\alpha,\beta)}\right)
u_{-2}^{k+1}
   +\left(\frac{1}{560}+\frac{1}{90}g_{1}^{(\alpha,\beta)}\right)
   u_{-2}^{k}
   +\frac{1}{90}\sum\limits_{\ell=2}^{k+1}
g_{\ell}^{(\alpha,\beta)}u_{-2}^{k+1-\ell}\\+
\left(\frac{3}{280}-\frac{3}{20}g_{0}^{(\alpha,\beta)}
\right)u_{-1}^{k+1}
-\left(\frac{3}{280}+\frac{3}{20}g_{1}^{(\alpha,\beta)}\right)u_{-1}^{k}
   -\frac{3}{20}\sum\limits_{\ell=2}^{k+1}
g_{\ell}^{(\alpha,\beta)}u_{-1}^{k+1-\ell}
\\-\left(\frac{1}{112}-\frac{3}{2}g_{0}^{(\alpha,\beta)}\right)u_{0}^{k+1}
   +\left(\frac{1}{112}+\frac{3}{2}g_{1}^{(\alpha,\beta)}\right)
   u_{0}^{k}
   +\frac{3}{2}\sum\limits_{\ell=2}^{k+1}
g_{\ell}^{(\alpha,\beta)}u_{0}^{k+1-\ell}\\
+\frac{1}{560}\tau f_{-2}^{k+\frac{1}{2}}
 -\frac{3}{280}\tau f_{-1}^{k+\frac{1}{2}}+\frac{1}{112}\tau
f_{0}^{k+\frac{1}{2}},
\vspace{0.1 cm}\\
-\left(\frac{1}{560}-\frac{1}{90}g_{0}^{(\alpha,\beta)}\right)u_{-1}^{k+1}
   +\left(\frac{1}{560}+\frac{1}{90}g_{1}^{(\alpha,\beta)}\right)
   u_{-1}^{k}
   +\frac{1}{90}\sum\limits_{\ell=2}^{k+1}
g_{\ell}^{(\alpha,\beta)}u_{-1}^{k+1-\ell}\\+
\left(\frac{3}{280}-\frac{3}{20}g_{0}^{(\alpha,\beta)}\right)u_{0}^{k+1}
-\left(\frac{3}{280}+\frac{3}{20}g_{1}^{(\alpha,\beta)}\right)u_{0}^{k}
   -\frac{3}{20}\sum\limits_{\ell=2}^{k+1}
g_{\ell}^{(\alpha,\beta)}u_{0}^{k+1-\ell}
\\+\frac{1}{560}\tau f_{-1}^{k+\frac{1}{2}}
 -\frac{3}{280}\tau f_{0}^{k+\frac{1}{2}},
\vspace{0.1 cm}\\
-\left(\frac{1}{560}-\frac{1}{90}g_{0}^{(\alpha,\beta)}\right)u_{0}^{k+1}
   +\left(\frac{1}{560}+\frac{1}{90}g_{1}^{(\alpha,\beta)}\right)
   u_{0}^{k} +\frac{1}{90}\sum\limits_{\ell=2}^{k+1}
g_{\ell}^{(\alpha,\beta)}u_{0}^{k+1-\ell}+
\frac{1}{560}\tau f_{0}^{k+\frac{1}{2}},
\vspace{0.2 cm}\\0
\vspace{0.2 cm}\\
      \vdots\vspace{0.2 cm}\\
       0 \vspace{0.1 cm} \\
-\left(\frac{1}{560}-\frac{1}{90}g_{0}^{(\alpha,\beta)}\right)u_{M}^{k+1}
   +\left(\frac{1}{560}+\frac{1}{90}g_{1}^{(\alpha,\beta)}\right)
   u_{M}^{k} +\frac{1}{90}\sum\limits_{\ell=2}^{k+1}
g_{\ell}^{(\alpha,\beta)}u_{M}^{k+1-\ell}+
\frac{1}{560}\tau f_{M}^{k+\frac{1}{2}},
\vspace{0.1 cm}\\
-\left(\frac{1}{560}-\frac{1}{90}g_{0}^{(\alpha,\beta)}\right)u_{M+1}^{k+1}
   +\left(\frac{1}{560}+\frac{1}{90}g_{1}^{(\alpha,\beta)}\right)
   u_{M+1}^{k}
   +\frac{1}{90}\sum\limits_{\ell=2}^{k+1}
g_{\ell}^{(\alpha,\beta)}u_{M+1}^{k+1-\ell}\\+
\left(\frac{3}{280}-\frac{3}{20}g_{0}^{(\alpha,\beta)}\right)u_{M}^{k+1}
-\left(\frac{3}{280}+\frac{3}{20}g_{1}^{(\alpha,\beta)}\right)u_{M}^{k}
   -\frac{3}{20}\sum\limits_{\ell=2}^{k+1}
g_{\ell}^{(\alpha,\beta)}u_{M}^{k+1-\ell}
\\+\frac{1}{560}\tau f_{M+1}^{k+\frac{1}{2}}
 -\frac{3}{280}\tau f_{M}^{k+\frac{1}{2}},\vspace{0.1 cm}\\
-\left(\frac{1}{560}-\frac{1}{90}g_{0}^{(\alpha,\beta)}\right)
u_{M+2}^{k+1}
   +\left(\frac{1}{560}+\frac{1}{90}g_{1}^{(\alpha,\beta)}\right)
   u_{M+2}^{k}
   +\frac{1}{90}\sum\limits_{\ell=2}^{k+1}
g_{\ell}^{(\alpha,\beta)}u_{M+2}^{k+1-\ell}\\+
\left(\frac{3}{280}-\frac{3}{20}g_{0}^{(\alpha,\beta)}
\right)u_{M+1}^{k+1}
-\left(\frac{3}{280}+\frac{3}{20}g_{1}^{(\alpha,\beta)}\right)u_{M+1}^{k}
   -\frac{3}{20}\sum\limits_{\ell=2}^{k+1}
g_{\ell}^{(\alpha,\beta)}u_{M+1}^{k+1-\ell}
\\-\left(\frac{1}{112}-\frac{3}{2}g_{0}^{(\alpha,\beta)}\right)u_{M}^{k+1}
   +\left(\frac{1}{112}+\frac{3}{2}g_{1}^{(\alpha,\beta)}\right)
   u_{M}^{k}
   +\frac{3}{2}\sum\limits_{\ell=2}^{k+1}
g_{\ell}^{(\alpha,\beta)}u_{M}^{k+1-\ell}\\
+\frac{1}{560}\tau f_{M+2}^{k+\frac{1}{2}}
 -\frac{3}{280}\tau f_{M+1}^{k+\frac{1}{2}}+\frac{1}{112}\tau
f_{M}^{k+\frac{1}{2}}
     \end{array}
    \right).
\end{array}
$$

{\it\bf Appendix II:}

(i) In the finite difference scheme (13), we use the following sixth-order extrapolation formulas for the ghost-point values:
$$
\begin{array}{lll} \displaystyle
 u_{-1}^k=6u_{0}^k-15u_{1}^k+20u_{2}^k-15u_{3}^k+6u_{4}^k-u_{5}^k+\mathcal {O}(h^6),
\end{array}
$$
and
$$
\begin{array}{lll} \displaystyle
 u_{M+1}^k=6u_{M}^k-15u_{M-1}^k+20u_{M-2}^k-15u_{M-3}^k+6u_{M-4}^k-u_{M-5}^k+\mathcal {O}(h^6).
\end{array}
$$

(ii) In the finite difference scheme (14), we use the following eighth-order extrapolation formulas for the ghost-point values:
$$
\begin{array}{lll} \displaystyle
 u_{-1}^k=8u_{0}^k-28u_{1}^k+56u_{2}^k-70u_{3}^k+56u_{4}^k-28u_{5}^k+8u_{6}^k-u_{7}^k+\mathcal {O}(h^8),
\end{array}
$$
$$
\begin{array}{lll} \displaystyle
 u_{-2}^k=36u_{0}^k-168u_{1}^k+378u_{2}^k-504u_{3}^k+420u_{4}^k-216u_{5}^k+63u_{6}^k-8u_{7}^k+\mathcal {O}(h^8),
\end{array}
$$
$$
\begin{array}{lll} \displaystyle
 u_{M+1}^k=8u_{M}^k-28u_{M-1}^k+56u_{M-2}^k-70u_{M-3}^k+56u_{M-4}^k-28u_{M-5}^k+8u_{M-6}^k-u_{M-7}^k+\mathcal {O}(h^8),
\end{array}
$$
and
$$
\begin{array}{lll} \displaystyle
 u_{M+2}^k=&\displaystyle36u_{M}^k-168u_{M-1}^k+378u_{M-2}^k-504u_{M-3}^k+420u_{M-4}^k-216u_{M-5}^k+63u_{M-6}^k
 \vspace{0.2 cm}\\&\displaystyle-8u_{M-7}^k+\mathcal {O}(h^8).
\end{array}
$$

  \end{document}